\documentclass[aop,MSNbibl,seceqn,citesort,dvips]{arximspdf}

\doi{10.1214/10-AOP612}
\volume{40}
\issue{1}
\pubyear{2012}
\firstpage{245}
\lastpage{279}

\makeatletter

\newcommand{\dd}{\mathrm{d}}

\newcommand{\xrightarrow}{\mathop{\hbox to 1cm{\rightarrowfill}}_}

\newcommand{\iint}{\int\!\!\int}

\newtheorem{theorem}{Theorem}
\newtheorem{lemma}{Lemma}
\newtheorem{proposition}{Proposition}
\newtheorem{corollary}{Corollary}

\newproclaim{remark}{Remark}

\newcommand{\er}{\mathbf{E}}
\newcommand{\pr}{\mathbf{P}}
\newcommand{\re}{\mathbb{R}}
\newcommand{\p}{\mathbb P}
\newcommand{\e}{\mathbb E}

\newcommand{\R}{\mathbb R}

\makeatother

\begin{document}
\begin{frontmatter}

\title{Fluctuation theory and exit systems for positive self-similar
Markov processes\thanksref{T1}}
\runtitle{Fluctuation theory and exit systems for pssMp}

\thankstext{T1}{Supported by EPSRC Grant EP/D045460/1, Royal Society
Grant RE-MA1004 and
ECOS-CONACYT Research Project M07-M01.}

\begin{aug}
\author[A]{\fnms{Lo\"{i}c} \snm{Chaumont}\ead[label=e1]{loic.chaumont@univ-angers.fr}},
\author[B]{\fnms{Andreas} \snm{Kyprianou}\ead[label=e2]{a.kyprianou@bath.ac.uk}},
\author[C]{\fnms{Juan Carlos} \snm{Pardo}\ead[label=e3]{jcpardo@cimat.mx}}\\ and
\author[C]{\fnms{V\'{i}ctor} \snm{Rivero}\corref{}\ead[label=e4]{rivero@cimat.mx}}

\runauthor{Chaumont, Kyprianou, Pardo and Rivero}
\affiliation{Universit\'{e} d'Angers, University of
Bath, CIMAT A.C. and CIMAT A.C.}
\address[A]{L. Chaumont\\
LAREMA, Universit\'{e} d'Angers\\
2 bd Lavoisier---49045 Angers Cedex 01\\
France\\
\printead{e1}} 
\address[B]{A. Kyprianou\\
Department of Mathematical Sciences\\
University of Bath\\
Bath, BA2 7AY\\
United Kingdom\\
\printead{e2}}
\address[C]{J. C. Pardo\\
V. Rivero\\
CIMAT A.C.\\
Calle Jalisco s/n\\
C.P.36240, Guanajuato, Gto.\\
Mexico\\
\printead{e3}\\
\hphantom{E-mail: }\printead*{e4}}
\end{aug}

\received{\smonth{12} \syear{2008}}

%
\begin{abstract}
For a positive self-similar Markov process, $X$, we construct a local
time for the random set, $\Theta$, of times where the process reaches
its past supremum. Using this local time we describe an exit system for
the excursions of $X$ out of its past supremum. Next, we define and
study the \textit{ladder process} $(R,H)$ associated to a positive
self-similar Markov process $X$, namely a bivariate Markov process with
a scaling property whose coordinates are the right inverse of the local
time of the random set $\Theta$ and the process $X$ sampled on the
local time scale. The process $(R,H)$ is described in terms of a ladder
process linked to the L\'{e}vy process associated to $X$ via
Lamperti's transformation. In the case where $X$ never hits $0$, and
the upward ladder height process is not arithmetic and has finite mean,
we prove the finite-dimensional convergence of $(R,H)$ as the starting
point of $X$ tends to $0$. Finally, we use these results to provide an
alternative proof to the weak convergence of $X$ as the starting point
tends to $0$. Our approach allows us to address two issues that
remained open in Caballero and Chaumont [\textit{Ann. Probab.}
\textbf{34} (2006) 1012--1034], namely, how to remove a~redundant
hypothesis and how
to provide a formula for the entrance law of $X$ in the case where the
underlying L\'{e}vy process oscillates.
\end{abstract}

%

\pdfkeywords{60G18, 60G17, 60J55, Entrance laws, exit systems,
excursion theory,
ladder processes, Lamperti's transformation, Levy processes,
self-similar Markov
processes}

\end{frontmatter}

\section{Introduction}
In recent years there has been a growing interest in the theory of
positive self-similar Markov processes (pssMp). Recall that a pssMp
$X=(X_t , t\geq0)$ is a right-continuous left-limited positive-valued
strong Markov process with the following scaling property. There exists
an $\alpha\neq0$ such that for any
$0<c<\infty$,
\[
\{ (cX_{tc^{-1/\alpha}}, t\geq
0),\p_x\}\stackrel{(\mathrm{d})}{=}\{ ( X_{t}, t\geq
0),\p_{cx}\},\qquad x>0.
\]
This class of processes has been introduced by Lamperti~\cite{La} in a
seminal paper where, among other interesting results, he established a
one-to-one relation between pssMp killed at $0$ and real-valued L\'
{e}vy processes. (Here we allow in the definition of a L\'{e}vy process
the additional possibility of being sent to a cemetery state after an
independent and exponentially distributed time).

Lamperti proved that any pssMp killed at $0$ is the exponential of a L\'
{e}vy process time changed by the right-continuous inverse of an
additive functional. We will refer to this relation as Lamperti's
transformation, and we will describe it in more detail in Section \ref
{preliminaries}. This relation allows one to embed the theory of L\'
{e}vy processes into that of pssMp. This embedding has proved to be a
powerful tool in unravelling the asymptotic behavior of pssMp (e.g.,~\cite{CP}),
in establishing various interesting identities (e.g.,~\cite{CKP})
and in linking these processes with other areas of applied
probability such as mathematical finance~\cite{yorbook}, continuous
state branching processes~\cite{KPbr} and fragmentation theory
\cite{Bfrag}, to name but a few.

Nevertheless implementing the theory of L\'{e}vy processes for such
purposes has never been as simple as one might hope
for on account of the time change, which relates
both classes together and destroys many of the convenient homogeneities
that are to be found in the theory of L\'{e}vy processes.
For example, it is known that a nondecreasing L\'{e}vy process with
finite mean grows linearly, owing to the law of large numbers, while a
pssMp associated via Lamperti's transformation to such a L\'{e}vy
process, growths with a polynomial rate whose order is given by the
index of self-similarity (see, e.g.,~\cite{beC} and~\cite{r2003}).

Our main objective in this paper is to shed light on fine properties
for the paths of pssMp, namely to establish a fluctuation theory, built
from the fluctuation theory of L\'{e}vy processes as well as classical
excursion theory (see, e.g.,~\cite{Be} and~\cite{Kbook} for
background). We will provide several new identities for pssMp and
present an alternative approach to that proposed by~\cite{beC,BeY}
and~\cite{CCh} for the existence of entrance laws for pssMp. The
latter allows us to address the problem of establishing an identity for
the entrance law at $0+$ for pssMp associated via Lamperti's
transformation to an oscillating L\'{e}vy process. To present our
results in more detail we need to introduce further notation and
preliminary results.


\section{Preliminaries and main results} \label{preliminaries}
Let $\mathbb{D}$ be the space of c\`{a}dl\`{a}g paths defined on
$[0,\infty)$, with values in $\R\cup\Delta$, where $\Delta$ is a
cemetery state. Each path $\omega\in\mathbb{D}$ is such that
$\omega_t=\Delta$, for any $t\ge\inf\{s\geq
0\dvtx\omega_s=\Delta\}:=\zeta(\omega)$. As usual we extend any function
$f\dvtx\re\to\re$ to $\re\cup\Delta$ by $f(\Delta)=0$. For each
Borel set $A$ we also define $\varsigma_A = \inf\{s>0\dvtx\omega_s \in
A\}$, writing, in particular, for convenience $\varsigma_0$ instead of
$\varsigma_{\{0\}}$. The space
$\mathbb{D}$ is endowed with the Skohorod topology and its Borel
$\sigma$-field. We will denote by $X$ the canonical process of the
coordinates.
Moreover, let $\pr$ be a reference probability
measure on~$\mathbb{D}$ under which the process $\xi$ is a L\'{e}vy
process; we will denote by $(\mathcal{G}_t, t\geq0)$ the complete
filtration generated by $\xi$. We will \textit{assume that the L\'{e}vy
process $(\xi,\pr)$ has an infinite lifetime}. 
Although this assumption remains in place throughout the paper, it
is not necessary and we indicate below how it may be removed.

Fix $\alpha\in\re\setminus\{0\}$, and let $(\p_x,x>0)$ be the laws
of the
$1/\alpha$-pssMp associated to $(\xi,\pr)$ via the Lamperti
representation. Formally, define
%
%
\begin{equation}\label{A}
A_t=\int^t_0\exp\{\alpha\xi_s\}\,\dd s,\qquad t\geq0,
\end{equation}
and let
$\tau(t)$ be its inverse,
\[
\tau(t)=\inf\{s>0\dvtx A_s>t\}
\]
with the
usual convention, $\inf\{\varnothing\}=\infty$. For $x>0$, we denote
by $\p_x$ the law of the process
\[
x\exp\bigl\{\xi_{\tau(tx^{-\alpha})}\bigr\},\qquad t>0.
\]
The Lamperti representation ensures that the laws $(\p_x, x>0)$ are
those of a~pssMp in the filtration $\{\mathcal{F}_{t}:=\mathcal
{G}_{\tau(t)}, t\geq0\}$ where index of self-similarity $1/\alpha$.
It follows that
$T_0=\inf\{t>0\dvtx X_{t}=0\}$ has the same law under $\p_x$
as~$x^{\alpha}A_{\infty}$ under $\pr$ with
\[
A_{\infty}=\int^{\infty}_0\exp\{\alpha\xi_s\}\,\dd s.
\]
Our assumption that $(\xi,\pr)$ has infinite lifetime implies that
the random variable
$A_{\infty}$ is finite a.s. or infinite a.s. according as $\lim_{t\to
\infty}\alpha\xi_t=-\infty$, a.s. or $\limsup_{t\to\infty}\alpha
\xi_t=\infty$, a.s. (see, e.g., Theorem 1 in~\cite{bysurvey}).
As a consequence, either $(X,\p_x)$ never hits $0$ a.s. or
continuously hits $0$ in a finite time a.s., independently of the
starting point $x>0$. Specifically in the latter case, the process
$(X,\p_x)$ does not jump to $0$. Lamperti~\cite{La} proved that all
pssMp that do not jump to $0$ can be constructed this way.

The assumption that $(\xi,\pr)$ has infinite lifetime can be removed
by killing the L\'{e}vy process $\xi$ at an independent and
exponentially distributed time with some parameter $q>0$ and then
applying the above explained transformation (Lamperti's transformation)
to the resulting killed L\'{e}vy process. Equivalently, one may kill
the pssMp $(X,\p)$, associated via Lamperti's transformation to a L\'
{e}vy process $(\xi,\pr)$, with infinite lifetime, by means of the
multiplicative functional
\[
\exp\biggl\{-q\int^t_{0}X^{-\alpha}_{s}\,\dd s\biggr\},\qquad
t\geq0.
\]
Using the Feymman--Kac formula to describe the infinitesimal generator
of the latter process it is readily seen that both procedures lead to
equivalent processes. Therefore,\vadjust{\goodbreak} we do not lose generality by making
the assumption that $(\xi,\pr)$ has infinite lifetime, as all our
results concern local properties of the process $(X,\p)$ or only make
sense for pssMp that never hit $0$.

Our main purpose is to study the paths of a pssMp by decomposing them
into the instants where it reaches its past-supremum. To this end, let
$(M_{t}, t\geq0)$ be the past supremum of $X$, $M_{t}=\sup\{X_{s},
0\leq s\leq t\}$, for $0\leq t\leq T_{0}$, and define the process $X$
\textit{reflected at its past supremum},
\[
\frac{M}{X}:=\biggl(\frac{M_{t}}{X_{t}}, 0\leq
t<T_{0}\biggr).
\]
It is easy to verify, using either the scaling and Markov properties or
Lamperti's transformation, that for $t,s\geq0$,
\[
M_{t+s}=X_{t} \biggl(\frac{M_{t}}{X_{t}}\vee\widetilde
{M}_{sX^{-\alpha}_{t}}\biggr),\qquad X_{t+s}=X_{t}\widetilde
{X}_{sX^{-\alpha}_{t}},
\]
where $\widetilde{M}$ (resp., $\widetilde{X}$) is a copy of $M$ (resp.,
of $X$) which is independent of $\mathcal{F}_{t}$. It follows that the
process $X$ reflected at its past supremum is not Markovian.
Nevertheless, using standard arguments it is easily established that
the process $Z$ defined by
\[
Z_{t}=\biggl(\frac{
M_{t}}{X_{t}}, M_{t}\biggr),\qquad0\leq t<T_{0},
\]
is a strong Markov process. Hence the random set of times $\Theta$
\[
\Theta=\biggl\{0\leq t< T_{0}\dvtx Z_{t}=\biggl(\frac
{M_{t}}{X_{t}},M_{t}\biggr)\in\{1\}\times\mathbb{R}_+\biggr\}
\]
is a homogeneous random set of $X$ in the sense of~\cite{maisonneuve},
but it is not regenerative in general because of its dependence on the
values of $M$. In~\cite{getoor} it has been proved, as a consequence
of the main result therein, that in a~rather general framework $(X,M)$
is a strong Markov process, and several functionals related to it have
been studied.

Our first aim is to describe the process $X$ at the instants of time in
$\Theta$ by means of the introduction of a \textit{local time}. To do so
we observe that because of Lamperti's transformation any element in
$\Theta$ is the image under the time change $\tau$ of some instant
where the underlying L\'{e}vy process reaches its past supremum. This
suggests that the random set $\Theta$ can be described by means of the
local time at $0$ of the underlying L\'{e}vy process reflected at its
past supremum. To make precise this idea we need to introduce further
notions related to the fluctuation theory of L\'{e}vy processes.

We recall that the process $\xi$ reflected at its past supremum,
$\overline{\xi}-\xi=(\sup_{s\leq t}\xi_{s}-\xi_{t}, t\geq0)$, is
a strong Markov process. Although all our results are true in general,
for brevity we will hereafter assume that:\vspace*{8pt}

\textit{$0$ is regular for $\overline{\xi}-\xi$ or}, \textit{equivalently}, \textit{that
$0$ is regular for $[0,\infty)$ for the
process~$\xi$.}\vadjust{\goodbreak}

Under these assumptions it is known that there is a continuous local
time at $0$ for $\overline{\xi}-\xi$, that as usual we will
denote by $L=(L_{s}, s\geq0)$ (see, e.g., Chapter IV in~\cite{Be}).
The instants where $\xi$ reaches its past supremum and the position of
$\xi$ at such times is described by
the so-called upward ladder time and height processes for $\xi$,
$(L^{-1}_{t}, \mbox{h}_{t}), t \geq0$, which are, respectively, defined
by
\[
L^{-1}_{t}=\cases{\inf\{s>0\dvtx L_{s}>t\}, &\quad$t<L_{\infty}$,\cr
\infty, &\quad$t\geq L_{\infty}$,}
\quad\mbox{and}\quad\mbox{h}_{t}=\cases{\xi_{L^{-1}_{t}}, &\quad
$t<L_{\infty}$,\cr
\infty, &\quad$t\geq L_{\infty}$.}
\]
It is known that the upward ladder process $(L^{-1},\mathrm{h})$ is a
bivariate L\'{e}vy process with increasing coordinates, and its
characteristics can be described in terms of $(\xi,\pr)$ (see, e.g.,
\cite{Be}, Chapter VI).
The downward ladder time and height processes are defined analogously,
replacing $\xi$ by its dual $(\widehat{\xi},\widehat{\pr}):=(-\xi
,\pr)$. We will assume that the downward and upward ladder time
subordinators are normalized such that their respective Laplace
exponents $\phi, \widehat{\phi}$, satisfy $\phi(1)=1=\widehat{\phi}(1)$.

We recall that there exists a constant $a\geq0$ such that the inverse
of the local time is given by
%
%
\begin{equation}\label{def:a}
L^{-1}_{t}=at+\sum_{s\leq t}\Delta L^{-1}_{s}, \qquad t\geq0,
\end{equation}
and
\[
\int^t_{0}\mathbf{1}_{\{\overline{\xi}_{s}-\xi_{s}=0\}}\,
\dd s=aL_{t},\qquad t\geq0.
\]
It is known that $a>0$ if and only if $0$ is irregular for $(-\infty
,0)$ for $\xi$. In that case the downward ladder time and height
processes are compound Poisson processes with inter-arrival rate $1/a$.
To see this recall, from the many statements that make up the
Wiener--Hopf factorization, that
\[
\phi(\lambda)\widehat{\phi}(\lambda)=\lambda,\qquad\lambda\geq0
\]
(see, e.g.,~\cite{Be}, Section VI.2). It follows that
\[
a=\lim_{\lambda\to\infty}\frac{\phi(\lambda)}{\lambda}=\lim
_{\lambda\to\infty}\frac{1}{\widehat{\phi}(\lambda)}=\frac
{1}{\Pi_{\hat{L}^{-1}}(0,\infty)},
\]
where $\Pi_{\hat{L}^{-1}}$ denotes the L\'{e}vy measure of
$\widehat{L}^{-1}$. Given that the downward ladder height subordinator
$\widehat{\mathrm h}$ stays constant in the same intervals where
$\widehat
{L}^{-1}$ does, the claim about the inter-arrival rate of $\widehat
{\mathrm h}$ follows.

We denote by $\{\epsilon_{t}, t\geq0\}$ the process of excursions of
$\xi$ from $\overline{\xi}$, namely
\[
\epsilon_{t}(s)=
\cases{\xi_{L^{-1}_{t-}}-\xi_{L^{-1}_{t-}+s}, 0\leq s\leq
L^{-1}_{t}-L^{-1}_{t-}, &\quad if
$L^{-1}_{t}-L^{-1}_{t-}>0$,\cr
\Delta, &\quad if $L^{-1}_{t}-L^{-1}_{t-}=0$.}
\]
It is well known that this process forms a Poisson point process on the
space of real valued c\`{a}dl\`{a}g paths with lifetime $\zeta$ and
whose intensity measure will be denoted by $\overline{n}$. This
measure is the so-called excursion measure from~$0$ for $\overline{\xi
}-\xi$. We will denote by $\epsilon$ the coordinate process under
$\overline{n}$. Under $\overline{n}$ the coordinate process has the
strong Markov property and its semigroup is that of $(\widehat{\xi
},\widehat{\pr})$ killed at $\varsigma_{(-\infty,0]}$, the first
hitting time of $(-\infty,0]$. Recall that for $A\subseteq\re$, we
denote by $\varsigma_{A}$ the first hitting time of the set $A$ for
the L\'{e}vy process $\xi$. Next, let $\widehat{\mathcal{V}}$ be the
measure over $\mathbb{R}_+$ defined by
%
%
\begin{equation}\label{dualrenewal}
\widehat{\mathcal{V}}(\dd y)=a\delta_{0}(\dd y)+\overline
{n}\biggl(\int^{\zeta}_{0}\mathbf{1}_{\{\epsilon(s)\in\dd y\}
}\,\dd s\biggr),\qquad y\geq0,
\end{equation}
and\vspace*{1pt} $\widehat{\mathcal{V}}(x)=a+\overline{n}(\int^{\zeta
}_{0}\mathbf{1}_{\{\epsilon(s)<x\}}\,\dd s)$ for $x>0$.
We recall that $\widehat{\mathcal{V}}(\dd x)$ equals the
potential measure $\er(\int^\infty_{0}\mathbf{1}_{\{\widehat
{\mathrm h}_{s}\in\dd x\}}\,\dd s)$, of the downward
ladder height subordinator $\widehat{\mathrm h}$. This is due to the fact
that
\[
\er\biggl(\int^\infty_{0}\mathbf{1}_{\{\widehat{\mathrm h}_{s}\in
\dd x\}}\,\dd s\biggr)=\overline{n}\biggl(\int^\zeta
_{0}\mathbf{1}_{\{\epsilon(s)\in\dd x\}}\,\dd s\biggr)\qquad
\mbox{on } x>0
\]
(see, e.g.,~\cite{Be}, exercise VI.6.5). Moreover, the potential
measure\break $\er(\int^\infty_{0}\mathbf{1}_{\{\widehat{\mathrm h}_{s}\in\dd
x\}}\,\dd s)$ has an atom at $0$ if
and only if $\widehat{\mathrm h}$ is a compound Poisson process. The
latter happens if and only if $0$ is irregular for $(-\infty,0)$ and
regular for $(0,\infty)$, for $\xi$, and then, in this case, the
aforementioned atom is of size $a$.

We introduce a new process
%
%
\begin{equation}\label{def:Y}\qquad
Y_{t}:=\cases{\displaystyle at+\sum_{s\leq t}\int
^{L^{-1}_{s}}_{L^{-1}_{s-}}\exp\{-\alpha(\xi_{L^{-1}_{s-}}-\xi
_{u})\}\,\dd s,
&\quad$0\leq t< L_{\infty}$,\cr
\infty, &\quad$t\geq L_{\infty}$.}
\end{equation}
It will be proved in Lemma~\ref{keylemma} that the process $Y$ is well
defined. Using the fact that the process $(L^{-1}, h, Y)$ is defined in
terms of functionals of the excursions from $0$ of the process $\xi$
reflected at it past supremum, standard arguments allow us to ensure
that the former process is a three-dimensional L\'{e}vy process whose
coordinates are subordinators, and they are not independent because
they have simultaneous jumps.

We will denote by $(\xi^{\uparrow}, \pr^{\uparrow})$ the process
obtained by conditioning $\xi$ to stay positive. We refer to \cite
{CD} and the references therein for further details on the construction
of $(\xi^{\uparrow}, \pr^{\uparrow})$. We will denote by $\pr
^{\dagger}$ the law of $\xi^{\uparrow}$ reflected in its future
infimum, $(\xi^{\uparrow}_{t}-\inf_{s\geq t}\xi^{\uparrow}_{s}$,
$t\geq0)$.

Given that the local time for $\overline{\xi}-\xi$ is an additive
functional in the filtration~$\mathcal{G}_{t}$, the process
\[
\widetilde{L}_{t}:=L_{\tau(t)}, \qquad t\geq0,
\]
is an additive functional in the filtration $(\mathcal{F}_{t}, t\geq
0)$ where $\mathcal{F}_{t}=\mathcal{G}_{\tau(t)}$ for $t\geq0$.
Moreover, $\widetilde{L}$ grows only on the instants where the process
$X$ reaches its past supremum, and thus a natural choice for the local
time of $\Theta$ would be $\widetilde{L}$. However, this additive
functional cannot be obtained as a limit of an occupation measure and
does not keep track of the values visited by the supremum.\vspace*{1pt}
It is for
these reasons that we will instead consider the process $L^{\Theta} =
(L^{\Theta}_t, t\geq0)$, defined by
%
%
\begin{equation}\label{def:LT}
L^{\Theta}_{t}=\int_{(0,t]}X^{\alpha}_{s}\,\dd\widetilde
{L}_{s},\qquad t\geq0.
\end{equation}
The latter is an additive functional of $Z = (M/X, M)$, and it is
carried by the set~$\Theta$. Unusually, however, the right-continuous
inverse of $L^{\Theta}$ is not a~subordinator on account of the fact
that $\Theta$ does not have the regenerative property as it depends on
the position of the past supremum. Nevertheless, $L^{\Theta}$ has
several properties common to local times of regenerative sets as is
demonstrated by the following result.
\begin{proposition}[(Local time)]\label{prop:localtime}
The constant $a$, defined in (\ref{def:a}), is such that
\[
\int^t_{0}\mathbf{1}_{\{M_{s}=X_{s}\}}\,\dd s=aL^{\Theta}_{t},\qquad
t\geq0.
\]
Assume that the underlying L\'{e}vy process $\xi$ is such that $0$ is
regular for $(0,\infty)$. The
function $\widehat{\mathcal{V}}$
normalizes the occupation measure in the following way:
\[
\lim_{\varepsilon\to0}\frac{1}{\widehat{\mathcal{V}}(\log
(1+\varepsilon))}\int^t_{0}\mathbf{1}_{\{
{M_{s}}/{X_{s}}\in[1,1+\varepsilon)\}}\,\dd s=L^{\Theta}_{t},
\]
uniformly over bounded intervals in $t$ where convergence is taken in
$\p_x$-probability, {for any starting point $x>0$}. We will therefore
refer to the process~$L^{\Theta}$ as the local time of the set $\Theta$.
\end{proposition}

An elementary but interesting remark is the following. By making
a~change of variables we obtain that under $\p_{1}$
\[
L^{\Theta}_{T_{0}}=\int^{L_{\infty}}_{0}e^{\alpha\mathrm{h}_{s}}\,\dd s,
\]
where we recall that $T_0=\inf\{t>0\dvtx X_t =0\}$.
Note that $L^\Theta_{T_0}$ is finite if and only if $L_{\infty
}<\infty$ or $L_{\infty}=\infty$ and $\alpha<0$. The terminal value
of the local time~$L^{\Theta}$ has the same law as the exponential
functional of a subordinator stopped at an independent exponential time
with some parameter $q\geq0$, where the value $q=0$ is allowed to
include the case where this random time is a.s. infinite.

In the following result we establish that there is a kernel associated
to $L^{\Theta}$, which we denote by $N^x$, with $x>0$ being the
initial value of the process of coordinates, so that $(L^{\Theta},
N^{\cdot})$ form an exit system for $\Theta$, which is an example of
the one introduced in~\cite{maisonneuve}.
\begin{theorem}\label{exitsystem}
Suppose that $F\dvtx\mathbb{R}_+^2\times\mathbb{D}\to\mathbb{R}_+$ is
a measurable function. Then for each $x>0$, define the kernel $N^x$ by
\[
N^{x}(F)=x^{-\alpha}\overline{n}\biggl(F\biggl(x, xe^{-\epsilon
({\zeta})}, xe^{\epsilon({\tau(t/x^{\alpha})})}, 0\leq t\leq
x^{\alpha}\int^\zeta_{0}e^{\alpha\epsilon(s)}\,\dd s
\biggr)\biggr).
\]
Let $\mathsf{G}$ be the set of left extrema of the excursion intervals
complementary to the random set $\Theta$, $D_{s}=\inf\{t>s\dvtx Z_t\in
\Theta\}$ for $s\geq0$, and $k_{\cdot}$ the killing operator. The
exit formula
\[
\e_{x}\biggl(\sum_{G\in\mathsf{G}}V_{G}F
\biggl(M_{G},M_{D_{G}},\biggl(\frac{M}{X}\circ\theta_{G}\biggr)\circ
k_{D_{G}}\biggr)\biggr)=\e_{x}\biggl(\int^\infty_{0}
\dd L^{\Theta}_{s}V_{s}N^{X_{s}}(F)\biggr)
\]
holds for every positive, left-continuous and $(\mathcal{F}_t, t\geq
0)$-adapted process $(V_t,\break t\geq0)$.
\end{theorem}
\begin{remark}
Observe that the kernel $N^{\cdot}(\cdot)$ has the following scaling
property: for $c>0$, the image of $N^{\cdot}(\cdot)$, under the
dilation $(cX_{tc^{-\alpha}}, t\geq0)$, equals $c^{\alpha}N^{c
\cdot}(\cdot)$.
\end{remark}

We are now ready to define the ascending ladder process associated to
the pssMp~$X$. The name arises from the analogous process appearing in
the fluctuation theory of L\'{e}vy processes.
\begin{theorem}[(Ascending ladder process)]\label{teo:upwardladderp}
Let $\{R_{t}, t\geq0\}$ be the right-continuous inverse of $L^{\Theta
}$, that is,
\[
R_{t}=\inf\{s>0\dvtx L^{\Theta}_{s}>t\},\qquad t\geq0,
\]
$H_{t}:=X_{R_{t}}$, $t\geq0$, and $K_{t}:=\int^{R_{t}}_{0}X^{-\alpha
}_{s}\,\dd s$, $t\geq0$. The following properties hold:
\begin{longlist}
\item The process $(K,R,H)$ has the same law under $\p_{x}$ as
the process
\[
\biggl\{\biggl(L^{-1}_{t}, x^{\alpha}\int_{(0,t]}e^{\alpha\mathrm{
h}_{s-}}\,\dd Y_{s}, xe^{\mathrm{h}_{t}}\biggr), t\geq0\biggr\},
\]
time changed by the inverse of the additive functional $(x^{\alpha
}\int^t_{0}e^{\alpha\mathrm{h}_{s}}\,\dd s, t\geq0)$, under~$\pr$.
\item The process $(R,H)$ is a Feller process in $[0,\infty
)\times(0,\infty)$ and has the following scaling property: for every
$c>0$, $((c^{\alpha}R_{tc^{-\alpha}},cH_{tc^{-\alpha
}}), t\geq0)$ issued from $(x_{1},x_{2})\in\re^2_{+}$
has the same law as $(R,H)$ issued from $(c^{\alpha}x_{1}, cx_{2})$.
\end{longlist}
\end{theorem}

We observe that Theorem~\ref{teo:upwardladderp} implies directly the
following corollary.\vadjust{\goodbreak}
\begin{corollary} Let $(R,H)$ be as in the previous theorem.
\begin{longlist}
\item The process $R$ is an increasing self-similar process
with index $1$. It is not a Markov process. Although, if $X$ has no
positive jumps, then $R$ has independent increments.
\item The process $H$ is the $1/\alpha$-increasing
self-similar Markov process which is obtained as the Lamperti transform
of the upward ladder height subordinator $\{{\mathrm h}_{t}, t\geq0\}$
associated to the L\'{e}vy process $\xi$.
\end{longlist}
\end{corollary}
\begin{remark}
It is important to mention that it is possible to state the analogues of
Proposition~\ref{prop:localtime} and Theorems~\ref{exitsystem} and
\ref{teo:upwardladderp} for the past infimum and~$X$ reflected in its
infumum, $\mathrm{I}_{t}=\inf_{0\leq s\leq t}X_{s}$, $t\geq0$,
$(\mathrm{I}_{t},X_{t}/\mathrm{I}_{t})$, $t\geq0$. These are easily
deduced from our results using the elementary fact that $X$ has the
same law as $1/\widehat{X}^{(-\alpha)}$, where $\widehat
{X}^{(-\alpha)}$ denotes the pssMp with self-similarity index
$1/(-\alpha)$, which is obtained by applying Lamperti's transformation
to $\widehat{\xi}:=-\xi$. We omit the details.
\end{remark}

Similar to fluctuation theory of L\'{e}vy processes, the process
$(R,H)$ is in general a simpler mathematical object to manipulate since
its coordinates are increasing processes and provide information about
$X$ at its running supremum.

Next, we will explain how the process $(R,H)$ can be used to provide an
alternative approach to those methods proposed by~\cite{beC,BeY}
and~\cite{CCh} with regard to establishing the existence of
entrance laws for pssMp. Working with the process $(R,H)$ has the
advantage of allowing us to give an explanation of the necessary and
sufficient conditions for the existence of entrance laws, and in doing
so we are able to remove an extra assumption in the main theorem of
\cite{CCh} as well as establish a general formula for the entrance law
at $0$ for $X$, thereby answering an open question from~\cite{CCh}.

Assume hereafter that $\alpha>0$. In~\cite{La} Lamperti remarked that
the Feller property at $0$ may fail for some pssMp and raised the
question of providing necessary and sufficient conditions for the
process $X$ to be a Feller process in $[0,\infty)$. For pssMp that hit
$0$ in a finite time this problem has been solved in complete
generality by Rivero~\cite{r2005,r2007} and Fitzsimmons
\cite{fitzsimmons}. It is known that for pssMp that never hit the state
$0$, the latter question is equivalent to studying the existence of
entrance laws. This task, as well as the proof of the weak convergence
of $(X, \p_x)$ as $x$ tends to zero, in the sense of finite-dimensional
distributions or in the sense of the weak convergence with respect to
the Skorohod topology, has been carried by Bertoin and Caballero in
\cite{beC}, Bertoin and Yor in~\cite{BeY} and Chaumont and Caballero in
\cite{CCh}. To explain their results and our results, we will assume
hereafter that the pssMp $X$ is such that
$\limsup_{t\to\infty}X_{t}=\infty$ a.s. which is equivalent to assuming
that the underlying L\'{e}vy process $\xi$ is such that
$\limsup_{t\to\infty}\xi_{t}=\infty$. In the case that the process
$\{{\mathrm h}_{s}, s\geq0\}$ is not arithmetic, they provide necessary
and sufficient conditions for the process $X$ to be Feller on
$[0,\infty)$ and to have weak convergence with respect to the
Skorohod\vspace*{1pt} topology as the starting point tends to~$0$.
These conditions are that $\er({\mathrm h}_{1})<\infty$ and $\er
(\log^+(\int^{\varsigma_{(1,\infty)}}_{0}e^{\alpha\xi _{s}}\,\dd
s))<\infty$. One of the main contributions of this paper is that we
will prove that in fact the sole condition $\er({\mathrm
h}_{1})<\infty$ is necessary and sufficient for the latter convergence
to hold.

Our first key observation in this direction is to remark that
Lamperti's question regarding the Feller nature of pssMp on $[0,\infty
)$ is equally applicable to the process $(R,H)$. The purpose of the
next theorem is to provide an answer to this question. The result can
be seen as an extension for the ladder process $(R,H)$ of the main
result in~\cite{beC}.
\begin{theorem}\label{limitladder}
If $\mathrm{h}$ is not arithmetic and $\mu_{+}=\er({\mathrm
h}_{1})<\infty$, then for every $t>0$ the bivariate measure $\p
_{x}(R_{t}\in\dd s, H_{t}\in\dd y)$
converges weakly as $x\to0+$ to a measure that we will denote by $\p
^{R,H}_{0+}(R_{t}\in\dd s, H_{t}\in\dd y)$ in $\re
^2_{+}$, which is such that for any measurable function $F\dvtx\mathbb
{R}_+^2\to\mathbb{R}_+$,
\[
\e^{R,H}_{0+}(F(R_{t}, H_{t}))=\frac{1}{\alpha\mu
_{+}}\er\biggl(F\biggl(\frac{t\widetilde{I}}{I_{\mathrm h}},\frac
{t^{1/\alpha}}{I^{1/\alpha}_{\mathrm h}}\biggr)\frac{1}{I_{\mathrm
h}}\biggr),\qquad t\geq0,
\]
where $I_{\mathrm h}=\int^\infty_{0}e^{-\alpha{\mathrm h}_{s}}\,\dd s$, and
$\widetilde{I}$ is the weak limit of $e^{-\alpha{\mathrm h}_{t}}\int
_{(0,t]}e^{\alpha{\mathrm h}_{s-}}\,\dd Y_{s}$, as \mbox{$t\to\infty$}.
Furthermore, $\widetilde{I}$ has the same law as $\int_0^{\infty}
\exp\{-\alpha\xi^{\uparrow}_s\}\,\dd s$. Finally, the process $\{
(R_{t},H_{t}), t\geq0\}$ converges in the sense of finite-dimensional
distributions as the starting point of $X$ tends to $0$.
\end{theorem}

It is implicit in Theorem~\ref{limitladder} that $\widetilde
{I}$ is a nondegenerate random variable.

In Corollary~\ref{th:exitformulares} we will represent the resolvent
of $X$ in terms of $(R,H)$ and then use Theorem~\ref{limitladder} to
prove the finite-dimensional convergence of the process $X$ as the
starting point tends to $0+$, and to obtain the formula for the
entrance law at $0+$ for $X$ in the following result. Weak convergence
with respect to the Skorohod topology will be then obtained by a
tightness argument.
\begin{theorem}\label{thm:mainconvresult}
Assume that $\xi$ is not arithmetic and that $\mu_{+}=\er({\mathrm
h}_{1})<\infty$. Then $\p_{x}$ converges weakly with respect to the
Skorohod topology, as the starting point $x$ tends to $0+$, toward a
probability measure $\p_{0+}$. The process $((X, \p_{x}), x\geq
0)$ is a strong Markov process and the one-dimensional law of~$X$
under $\p_{0+}$ is determined by
%
%
\begin{equation}\label{nentlaw}
\e_{0+}(f(X_{t}))=\int^\infty_{0}f\biggl(\frac
{t^{1/\alpha}}{x^{1/\alpha}}\biggr)\frac{1}{x} \eta(\dd x),
\end{equation}
where $\eta$ is the measure defined by
\[
\eta(f)=\frac{1}{\alpha\mu_{+}}\int_{\re^{3}_{+}} \pr(\widetilde
{I}\in\dd t)\widehat{\mathcal{V}}(\dd x) \pr^{\dagger
}_{x}\biggl(\int^{\varsigma_{0}}_{0}e^{-\alpha\xi_{u}}\,\dd u\in\dd
s\biggr)f\bigl(e^{\alpha x}(t+s
)\bigr)
\]
and $\int^{\infty}_{0}x^{-1}\eta(\dd x)=1$.
\end{theorem}

Roughly speaking, our approach to proving Theorem \ref
{thm:mainconvresult} relies on the idea that we need first to prove the
convergence of its ladder height process.
%
\begin{remark}
In~\cite{CCh}, it has been proved that if $\e(\log^+\int
_0^{\varsigma_{(1,\infty)}}\exp\xi_s \,\dd s)=\infty$, then
$\p_x$ converges weakly toward the degenerated process $X_\cdot
\equiv0$. Under conditions of Theorem~\ref{thm:mainconvresult}, the weak
limit is not degenerated; therefore, a simple argument by contradiction
shows that, within the context of Theorem~\ref{thm:mainconvresult}, it
necessarily holds that
$\e(\log^+\int_0^{\varsigma_{(1,\infty)}}\exp\xi_s\, \dd s)<\infty$. Note
also that the conditions of Theorem \ref
{thm:mainconvresult}, that is, that
$\xi$ is not arithmetic, and $\er({\mathrm h}_{1})<\infty$ are equivalent
to the following:
\[
\mbox{$\xi$ is not arithmetic}\quad\mbox{and}\quad
\cases{\mbox{either $0<\er(\xi_1)\leq\er(|\xi_1|)<\infty$,}\cr
\mbox{or $\er(|\xi_1|)<\infty$, $\er(\xi_1)=0$ and $J<\infty$,}}
\]
where $ J=\int_{[1,\infty)}\frac{x{\pi}(x,\infty)
\,\dd x}{1+\int_0^x\dd y \int_y^\infty
{\pi}(-\infty,-z) \,\dd z}$ and $\pi$ is the L\'{e}vy measure of $\xi$
(see~\cite{CCh}, Section 2.1).
\end{remark}
%
\begin{remark}
When $\mathbf{E}(\mathrm{h}_1)<\infty$ and $\xi$ drifts to $+\infty
$, that is, $\er({\xi}_{1})>0$, an expression for the entrance law
under $\p_{0+}$ has already been obtained by Bertoin and Yor
\cite{BeY}. More precisely, for every $t>0$ and measurable function
$f\dvtx\R
^+\to\R^+$, we have
%
%
\begin{equation}\label{entlaw}
\e_{0+}(f(X_t))=\frac{1}{\alpha\er(\xi_1)}\er
\bigl(I^{-1}f\bigl((t/I)^{1/\alpha}\bigr)\bigr).
\end{equation}
In this case and under some mild technical conditions on $\xi$ which
can be found in~\cite{Ch}, formula (\ref{nentlaw}) can be recovered
using pathwise arguments as follows: Let $m$ be the unique time when
$\xi$ reaches its overall minimum. It is well known (cf.~\cite{Mi,Ch})
that $(\xi_t, 0\leq t < m)$ and $(\xi_{t+m}-\xi_m, t \ge0)$
are independent, and the latter process has the same law as $\xi$
conditioned to stay positive. Moreover, from Lemma~8 of~\cite{CR}, the
law of the pre-minimum process is characterized by
\[
\er\bigl(H(\xi_{t}-\xi_m, 0\le t < m)\bigr)=\widehat{\kappa}\int
_0^{\infty} \widehat{\mathcal{V}}(\dd x)\pr_x^{\downarrow
}\bigl(H(\xi_{t}, 0\le t < \zeta)\bigr),
\]
where $\widehat{\kappa}$ is the killing rate of the subordinator
$\widehat{\mathrm h}$, and $\pr_x^{\downarrow}$ is the law of the L\'
{e}vy process
$\xi$ starting from $x$ and conditioned to hit 0 continuously. Then we write
\[
I=\int_0^{\infty}e^{-\alpha\xi_s}\,\dd s=e^{-\alpha\xi_m}
\biggl(\int_0^{m}e^{-\alpha(\xi_s-\xi_m)}\,\dd s+\int
_0^{\infty}e^{-\alpha(\xi_{m+s}-\xi_m)}\,\dd s\biggr).
\]
From this identity we obtain that for any measurable function $g\dvtx\R
^+\to\R^+$,
\[
\er(g(I))=\widehat{\kappa}\int_0^{\infty}\pr(I^{\uparrow
}\in\dd r)\int_0^{\infty} \widehat{\mathcal{V}}(\mathrm
{d}x)\pr_x^{\downarrow}\biggl(\int_0^{\zeta}e^{-\alpha\xi
_s}\,\dd s\in\dd u\biggr)g\bigl(e^{\alpha x}(r+u)\bigr),
\]
where $I^{\uparrow}=\int_0^{\infty}e^{-\alpha(\xi_{m+s}-\xi_
m)}\,\dd s\stackrel{(\mathrm{d})}{=}\int_0^{\infty}e^{-\alpha\xi
^\uparrow
_{s}}\,\dd s$. Taking $g(x)=(\alpha\er(\xi
_1))^{-1}\times f(x^{-1/\alpha})x^{-1}$, it gives
\begin{eqnarray*}
&&\frac{1}{\alpha\er(\xi_1)}\er\bigl(I^{-1}f\bigl((1/I)^{1/\alpha
}\bigr)\bigr)\\[-2pt]
&&\qquad=\frac{\widehat{\kappa}}{\alpha\er(\xi_1)}\int_0^{\infty}\pr
(I^{\uparrow}\in\dd r)\int_0^{\infty} \widehat{\mathcal
{V}}(\dd x)\pr_x^{\downarrow}\biggl(\int_0^{\zeta}e^{-\alpha
\xi_s}\,\dd s\in\dd u\biggr)\\[-2pt]
&&\qquad\quad{}\times f\bigl(e^{-x}(r+u)^{-1/\alpha
}\bigr)\frac{1}{e^{\alpha x}(r+u)}.
\end{eqnarray*}
It follows from Theorem 5 in~\cite{Ch} that the law of the canonical
process killed at its first hitting time of $0$, $\varsigma_0$, under
$\pr^{\dagger}_{x}$ equals that of the process $\xi$ issued from $x$
and conditioned to hit $0$ continuously under $\pr^{\downarrow}_{x}$.
Finally, we use the fact that $\er(\xi_{1})=\widehat{\kappa}\er
(\mathrm{h}_{1})$, which is a simple consequence of the Wiener--Hopf
factorization (see, e.g.,~\cite{doney}, Corollary 4.4) to recover
formula (\ref{nentlaw}) in Theorem~\ref{thm:mainconvresult}.\vspace*{-2pt}
\end{remark}

The rest of this paper is organized as follows. Section~\ref{ladders}
is mainly devoted to prove Proposition~\ref{prop:localtime} and
Theorem~\ref{exitsystem}. In proving Proposition~\ref{prop:localtime}
we will use that the local time at $0$ of a L\'{e}vy process reflected
in its past supremum can be approximated by an occupation time
functional. This result is of interest in itself and, to the best of
our knowledge this cannot be found in the literature, so we have
included a proof. Next we use the excursion theory for L\'{e}vy
processes reflected in its past supremum to establish Theorem \ref
{exitsystem}. The main purpose of Section~\ref{section:Pthm2} is to
prove Theorem~\ref{teo:upwardladderp}. To this end we will establish
an elementary but key lemma that allows us to describe the time change
appearing in Lamperti's transformation in terms of the excursions of
the process $\xi$ out of its past supremum. In Section \ref
{section:res} we describe the $q$-resolvent of $X$ in terms of the
ladder process $(R,H)$, which will be useful in the proof of the
convergence results. As an application of these results, in
Section~\ref{sect:appli} we will prove Theorems~\ref{limitladder} and
\ref{thm:mainconvresult}.\vspace*{-2pt}

\section{\texorpdfstring{Proofs of Proposition \protect\ref{prop:localtime} and Theorem \protect\ref{exitsystem}}
{Proofs of Proposition 1 and Theorem 1}}\label{ladders}

The proof of Proposition~\ref{prop:localtime} needs the following
analogous result for L\'{e}vy processes which is of interest in itself.\vspace*{-2pt}
\begin{lemma}\label{localtimelevy}
Let $\xi$ be a L\'{e}vy process for which $0$ is regular for
$(0,\infty)$. Then
\[
\lim_{\varepsilon\downarrow0}\frac{1}{\widehat{\mathcal
{V}}(\varepsilon)}\int^{L^{-1}_{t}}_{0}\mathbf{1}_{\{\overline{\xi
}_{s}-\xi_{s}<\varepsilon\}}\,\dd s=t\wedge L_{\infty},\vadjust{\goodbreak}
\]
uniformly over bounded intervals of $t$ in the $L^2$-norm. Furthermore,
\[
\lim_{\varepsilon\downarrow0}\frac{1}{\widehat{\mathcal
{V}}(\varepsilon)}\int^{t}_{0}\mathbf{1}_{\{\overline{\xi}_{s}-\xi
_{s}<\varepsilon\}}\,\dd s=L_{t},
\]
uniformly over bounded intervals in probability.
\end{lemma}
\begin{pf}
In the case where $\widehat{\mathcal{V}}(0)=a>0$, we know that
\[
\int^t_{0}\mathbf{1}_{\{\overline{\xi}_{s}-\xi_{s}=0\}}\,\dd
s=aL_{t},\qquad t>0.
\]
Using this fact it is readily seen that the claims are true, so we can
restrict ourselves to the case $a=0$. The proof of the first assertion
in this lemma follows the basic steps of the proof of this result for
L\'{e}vy processes with no negative jumps due to Duquesne and Le
Gall~\cite{DL}.

Let $\mathcal{N}$ be a Poisson random measure on $\mathbb{R}_+\times
\mathbb{D}(\mathbb{R}_+,\re)$ with intensity $\dd t\times
\overline{n}(\dd\epsilon)$. For every $t>0$, put
\[
J_{\varepsilon}(t):=\frac{1}{\widehat{\mathcal{V}}(\varepsilon
)}\int\mathcal{N}(\dd u, \dd\epsilon)\cdot\mathbf
{1}_{\{u\leq t\}}\int^\zeta_{0}\mathbf{1}_{\{[0,\varepsilon)\}
}(\epsilon(s))\,\dd s.
\]
It follows that
\[
\er(J_{\varepsilon}(t))=\frac{t}{\widehat{\mathcal{V}}(\varepsilon
)}\overline{n}\biggl(\int^\zeta_{0}\mathbf{1}_{\{[0,\varepsilon)\}
}(\epsilon(s))\,\dd s\biggr)=t,
\]
where the second equality follows from the definition of $\widehat
{\mathcal{V}}(\varepsilon)$. Furthermore,
\[
\er((J_{\varepsilon}(t))^2)=(\er
(J_{\varepsilon}(t)))^2+\frac{t}{(\widehat
{\mathcal{V}}(\varepsilon))^2}\overline{n}\biggl(\biggl(\int^\zeta
_{0}\mathbf{1}_{\{\epsilon(s)\in[0,\varepsilon)\}}\,\dd s\biggr)^2\biggr),
\]
the latter equality can be verified using that for each $\varepsilon
>0$, the process $(J_{\varepsilon}(t)$, $t\geq0)$ is a subordinator
whose L\'{e}vy measure is the image measure of $\overline{n}$ under
the mapping $\epsilon\mapsto\widehat{\mathcal{V}}(\varepsilon
)^{-1}\int^\zeta_{0}\mathbf{1}_{\{[0,\varepsilon)\}}(\epsilon
(s))\,\dd s$, and thus the second moment is obtained by
differentiating twice the Laplace transform and using the L\'
{e}vy--Khintchine formula. The right-most term in the latter identity
can be estimated as follows:
\begin{eqnarray*}
&&
\overline{n}\biggl(\biggl(\int^\zeta_{0}\mathbf{1}_{\{\epsilon
(s)\in[0,\varepsilon)\}}\,\dd s\biggr)^2\biggr)\\
&&\qquad=2\overline
{n}\biggl(\int_{0\leq s\leq t \leq\zeta}\mathbf{1}_{\{\epsilon
(s)\in[0,\varepsilon)\}}\mathbf{1}_{\{\epsilon(t)\in[0,\varepsilon
)\}}\,\dd s\,\dd t\biggr)\\
&&\qquad=2\overline{n}\biggl(\int^{\zeta}_{0}\mathbf{1}_{\{\epsilon
(s)\in
[0,\varepsilon)\}}\er_{-\epsilon(s)}\biggl(\int^{\varsigma
_{(0,\infty)}}_{0}\mathbf{1}_{\{\xi(t)\in(-\varepsilon,0]\}
}\,\dd t\biggr)\,\dd s\biggr)\\
&&\qquad\leq2\widehat{\mathcal{V}}(\varepsilon) \sup_{y\in
[0,\varepsilon
)}\er_{-y}\biggl(\int^{\varsigma_{(0,\infty)}}_{0}\mathbf{1}_{\{
\xi(t)\in(-\varepsilon,0]\}}\,\dd t\biggr),
\end{eqnarray*}
where we have used the Markov property for excursions in the second equality.

By Theorem VI.20 in~\cite{Be} it is known that there exists a
constant, say $k$, such that for $y<\varepsilon$
\begin{eqnarray*}
\er_{-y}\biggl(\int^{\varsigma_{(0,\infty)}}_{0}\mathbf{1}_{\{\xi
(t)\in(-\varepsilon,0]\}}\,\dd t\biggr)&=&k\int_{[0,y)}\mathcal
{V}(\dd z)\int_{[0,\infty)}\widehat{\mathcal{V}}(\mathrm
{d}x)\mathbf{1}_{(-\varepsilon,0]}(-y+z-x)\\
&=&k\int_{[0,y)}\mathcal{V}(\dd z)\int_{[0,\infty)}\widehat
{\mathcal{V}}(\dd x)\mathbf{1}_{[0,\varepsilon)}(x+y-z)\\
&\leq& k\mathcal{V}[0,y)\widehat{\mathcal{V}}[0,\varepsilon)\\
&\leq& k\widehat{\mathcal{V}}(\varepsilon)\mathcal{V}[0,\varepsilon
)\\
&=&o(\widehat{\mathcal{V}}(\varepsilon)),
\end{eqnarray*}
where the second inequality follows from the fact that $\widehat
{\mathcal{V}}(z)\leq\widehat{\mathcal{V}}(\varepsilon)$, $0<z\leq
\varepsilon$, and the fifth from the fact that $\mathcal
{V}[0,\epsilon)\to0$ as $\epsilon\to0+$, because $0$ is regular for
$(0,\infty)$ (which implies that the upward subordinator is not a
compound Poisson process and so that its renewal measure $\mathcal{V}$
does not have an atom at $0$). It therefore follows that
%
%
\begin{equation}\label{eq:ocmes}\overline{n}\biggl(\biggl(\int^\zeta
_{0}\mathbf{1}_{\{\epsilon(s)\in[0,\varepsilon)\}}\,\dd s\biggr)^2\biggr
)= o((\widehat{\mathcal
{V}}(\varepsilon))^2).
\end{equation}
These estimates allow us to ensure that
\[
\lim_{\varepsilon\to0}\er\bigl(\bigl(J_{\varepsilon}(t)-t
\bigr)^2\bigr)=0.
\]
Moreover, thanks to the fact that $(J_{\varepsilon}(t)-t, t\geq0)$ is
a martingale, we can apply Doob's inequality to deduce that
\[
\lim_{\varepsilon\to0}\er\Bigl(\sup_{s\leq t}\bigl(J_{\varepsilon
}(s)-s\bigr)^2\Bigr)=0.
\]
The first assertion of the theorem follows since the pair
\[
\biggl(\frac{1}{\widehat{\mathcal{V}}(\varepsilon)}\int
^{L^{-1}_{t}}_{0}\mathbf{1}_{\{\overline{\xi}_{s}-\xi_{s}<\epsilon
\}}\,\dd s, L_{\infty}\biggr)
\]
has the same law as $(J_{\varepsilon}(t\wedge\nu), \nu)$ where
\mbox{$\nu=\inf\{t>0\dvtx\mathcal{N}(0,t]\times\{\zeta=\infty\}\geq1\}$}.

Now to prove the second assertion we fix $t>0$, let $\varepsilon_{1},
\delta>0$ and take $T>0$ large enough such that $\p
(L_{t}>T)<\varepsilon_{1}/3$. It follows using the inequalities
$L^{-1}_{L_{s}}\geq s\geq L^{-1}_{L_{s}-}$ that
%
%
\begin{eqnarray}\label{eq:ocnvunifp}
&&\pr\biggl(\sup_{s\leq t}\biggl|\frac{1}{\widehat{\mathcal
{V}}(\varepsilon)}\int^s_{0}\mathbf{1}_{\{\overline{\xi}_{u}-\xi
_{u}<\varepsilon\}}\,\dd u-L_{s}\biggr|>\delta\biggr)
\nonumber\\
&&\qquad\leq\pr\biggl(\sup_{s\leq t}\biggl|\frac{1}{\widehat{\mathcal
{V}}(\varepsilon)}\int^{L^{-1}_{L_{s}}}_{0}\mathbf{1}_{\{\overline
{\xi}_{u}-\xi_{u}<\varepsilon\}}\,\dd u-L_{s}\biggr|>\delta/2,
L_{t}<T\biggr)\nonumber\\
&&\qquad\quad{} +\pr\biggl(\sup_{s\leq t}\frac{1}{\widehat{\mathcal
{V}}(\varepsilon)}\int^{L^{-1}_{L_{s}}}_{L^{-1}_{L_{s}-}}\mathbf
{1}_{\{\overline{\xi}_{u}-\xi_{u}<\varepsilon\}}\,\dd u>\delta/2,
L_{t}<T\biggr)+\varepsilon_{1}/3\\
&&\qquad\leq\pr\biggl(\sup_{s\leq T}\biggl|\frac{1}{\widehat{\mathcal
{V}}(\varepsilon)}\int^{L^{-1}_{s}}_{0}\mathbf{1}_{\{\overline{\xi
}_{u}-\xi_{u}<\varepsilon\}}\,\dd u-s\biggr|>\delta/2
\biggr)\nonumber\\
&&\qquad\quad{} +\pr\biggl(\sup_{s\leq T}\int
^{L^{-1}_{s}}_{L^{-1}_{s-}}\mathbf{1}_{\{\overline{\xi}_{u}-\xi
_{u}<\varepsilon\}}\,\dd u>\delta\widehat{\mathcal{V}}(\varepsilon
)/2\biggr)+\varepsilon_{1}/3.\nonumber
\end{eqnarray}
It follows from the first assertion in Lemma~\ref{localtimelevy} that
$\varepsilon$ can be chosen so that the first term in the right-hand
side in inequality (\ref{eq:ocnvunifp}) is smaller than $\varepsilon
_{1}/3$. Moreover, as the random objects
\[
\int^{L^{-1}_{s}}_{L^{-1}_{s-}}\mathbf{1}_{\{\overline{\xi}_{u}-\xi
_{u}<\varepsilon\}}\,\dd u,\qquad s\geq0,
\]
are the values of the points in a Poisson point process in $\mathbb
{R}_+$ whose intensity measure is the image of $\overline{n}$ under
the mapping $\epsilon\mapsto\int^{\zeta}_{0}\mathbf{1}_{\{\epsilon
(u)<\varepsilon\}}\,\dd u$, it follows that
%
%
\begin{eqnarray}\label{eq:extremal}
&&\pr\biggl(\sup_{s\leq T}\int
^{L^{-1}_{s}}_{L^{-1}_{s-}}\mathbf{1}_{\{\overline{\xi}_{u}-\xi
_{u}<\varepsilon\}}\,\dd u>\delta\widehat{\mathcal{V}}(\varepsilon
)/2\biggr)\nonumber\\[-8pt]\\[-8pt]
&&\qquad=1-\exp\biggl\{-T\overline{n}\biggl(\int^{\zeta
}_{0}\mathbf{1}_{\{\epsilon(u)<\varepsilon\}}\,\dd u>\delta
\widehat{\mathcal{V}}(\varepsilon)/2\biggr)\biggr\}.\nonumber
\end{eqnarray}
From the Markov inequality and (\ref{eq:ocmes}) we have also that
\begin{eqnarray*}
\overline{n}\biggl(\int^{\zeta}_{0}\mathbf{1}_{\{\epsilon
(u)<\varepsilon\}}\,\dd u>\delta\widehat{\mathcal{V}}(\varepsilon
)/2\biggr)&\leq&\frac{4}{\delta^2 \widehat{\mathcal{V}}(\varepsilon
)^2}\overline{n}\biggl(\biggl(\int^{\zeta}_{0}\mathbf{1}_{\{
\epsilon(u)<\varepsilon\}}\,\dd u\biggr)^2\biggr)\\
&=&o(1)\qquad\mbox
{as } \varepsilon\to0.
\end{eqnarray*}
It follows then that by taking $\varepsilon$ small enough the
right-hand term in (\ref{eq:extremal}) can be made smaller
than $\varepsilon_{1}/3$, which finishes the proof of the second claim.
\end{pf}
\begin{pf*}{Proof of Proposition~\ref{prop:localtime}}
To prove the first claim we recall that the constant $a$ is such that
\[
\int^t_{0}\mathbf{1}_{\{\overline{\xi}_{s}-\xi_{s}=0\}}\,\dd
s=aL_{t},\qquad
t\geq0.
\]
Hence, by making a time change we get that for every $t>0$,
%
%
\begin{eqnarray}\label{eq:idLT}
\int^t_{0}\!\mathbf{1}_{\{M_{s}=X_{s}\}}\,\dd s&=&x^{\alpha}\!\int
^{\tau(t/x^{\alpha})}_{0}\!\mathbf{1}_{\{\overline{\xi}_{u}-\xi
_{u}=0\}}e^{\alpha\xi_{u}}\,\dd u=ax^{\alpha}\!\int_{(0,\tau(t/x^{\alpha
})]}\!e^{\alpha\xi
_{u}}\,\dd L_{u}\nonumber\hspace*{-35pt}\\[-4pt]\\[-12pt]
&=&a\int_{(0,t]}\bigl(xe^{\xi_{\tau(s/x^{\alpha
})}}\bigr)^{\alpha}\,\dd L_{\tau(s/x^{\alpha})}=aL^{\Theta}_{t}
\nonumber
\end{eqnarray}
under $\p_{x}$.
Now to prove the second claim we observe first that, as in the proof of
Lemma~\ref{localtimelevy}, we can restrict ourselves to the case where $a=0$.

For notational convenience, and without loss of generality thanks to
the self-similarity of $X$, we will assume that $X$ is issued from $1$.
By applying Lamperti's transformation and making a time change we
obtain the elementary inequalities
%
%
\begin{eqnarray}
&&\frac{1}{\widehat{\mathcal{V}}(\log(1+\varepsilon))(1+\varepsilon
)^{\alpha}}\int^{\tau(t)}_{0}\mathbf{1}_{\{\overline{\xi
}_{s}-{\xi}_{s}\in[0,\log(1+\varepsilon))\}}e^{\alpha
\overline{\xi}_{s}}\,\dd s\nonumber\\
&&\qquad\leq\frac{1}{\widehat{\mathcal{V}}(\log(1+\varepsilon))}\int
^t_{0}\mathbf{1}_{\{{M_{s}}/{X_{s}}\in[1,1+\varepsilon
)\}}\,\dd s\nonumber\\[-8pt]\\[-8pt]
&&\qquad=\frac{1}{\widehat{\mathcal{V}}(\log(1+\varepsilon))}\int^{\tau
(t)}_{0}\mathbf{1}_{\{\overline{\xi}_{s}-{\xi}_{s}\in[0,\log
(1+\varepsilon))\}}e^{-\alpha(\overline{\xi}_{s}-\xi
_{s})}e^{\alpha\overline{\xi}_{s}}\,\dd s\nonumber\\
&&\qquad\leq\frac{1}{\widehat{\mathcal{V}}(\log(1+\varepsilon))}\int
^{\tau(t)}_{0}\mathbf{1}_{\{\overline{\xi}_{s}-{\xi}_{s}\in
[0,\log(1+\varepsilon))\}}e^{\alpha\overline{\xi
}_{s}}\,\dd s.\nonumber
\end{eqnarray}
Let $\delta,t>0$ fixed. We infer the following inequalities:
%
%
\begin{eqnarray}\quad
&&\p_{1}\biggl(\sup_{r\leq t}\biggl|\frac{1}{\widehat{\mathcal
{V}}(\log(1+\varepsilon))}\int^r_{0}\mathbf{1}_{\{
{M_{u}}/{X_{u}}\in[1,1+\varepsilon)\}}\,\dd u - L^{\Theta
}_{r}\biggr|>\delta\biggr)\nonumber\\
&&\qquad\leq\pr\biggl(\frac{1}{\widehat{\mathcal{V}}(\log(1+\varepsilon
))}\int^{\tau(r)}_{0}\mathbf{1}_{\{\overline{\xi}_{s}-{\xi
}_{s}\in[0,\log(1+\varepsilon))\}}e^{\alpha\overline{\xi
}_{s}}\,\dd s-L^{\Theta}_{r}>\delta,\nonumber\\
&&\hspace*{240pt} \mbox{for some } r\leq
t\biggr)\\
&&\qquad\quad{} +\pr\biggl(L^{\Theta}_{r}-\frac{1}{(1+\varepsilon
)^{\alpha
}\widehat{\mathcal{V}}(\log(1+\varepsilon))}\nonumber\\
&&\hspace*{84.7pt}{}\times\int^{\tau
(r)}_{0}\mathbf{1}_{\{\overline{\xi}_{s}-{\xi}_{s}\in[0,\log
(1+\varepsilon))\}}e^{\alpha\overline{\xi}_{s}}\,\dd
s>\delta\mbox{,
for some } r\leq t\biggr).
\nonumber
\end{eqnarray}
Next we will prove that the probability in the first term on the
right-hand side tends to $0$ as $\varepsilon\to0$. The arguments used
to prove that the second one tends to $0$ are similar, so we omit them.
Consider the event
\begin{eqnarray*}
\mathcal{A}&=&\biggl\{\sup_{s\leq T}\biggl|\frac{1}{\widehat{\mathcal
{V}}(\log(1+\varepsilon))}\int^{s}_{0}\mathbf{1}_{\{\overline
{\xi}_{u}-{\xi}_{u}\in[0,\log(1+\varepsilon))\}}\,\dd
u-L_{s}\biggr|<\delta_{1},\\
&&\hspace*{128.5pt} \tau(t)<T, (2e^{\alpha\overline{\xi
}_{T}}-1)<\delta/\delta_{1}\biggr\}
\end{eqnarray*}
for $\delta_{1},\varepsilon,T>0$. From Lemma~\ref{localtimelevy} and
standard arguments it follows that $\delta_{1}, T$ and $\varepsilon$
can be chosen so that the probability of the event $\mathcal{A}$ is
arbitrarily close to $1$. By integrating by parts twice, we have that
on $\mathcal{A}$
\begin{eqnarray*}
&&\frac{1}{\widehat{\mathcal{V}}(\log(1+\varepsilon))}\int^{\tau
(r)}_{0}\mathbf{1}_{\{\overline{\xi}_{s}-{\xi}_{s}\in[0,\log
(1+\varepsilon))\}}e^{\alpha\overline{\xi}_{s}}\,\dd s\\
&&\qquad= e^{\alpha\overline{\xi}_{\tau(r)}}\frac{1}{\widehat{\mathcal
{V}}(\log(1+\varepsilon))}\int^{\tau(r)}_{0}\mathbf{1}_{\{
\overline{\xi}_{s}-{\xi}_{s}\in[0,\log(1+\varepsilon))\}
}\,\dd s\\
&&\qquad\quad{} -\int_{(0,\tau(r)]}\biggl(\frac{1}{\widehat{\mathcal
{V}}(\log(1+\varepsilon))}\int^{s}_{0}\mathbf{1}_{\{\overline
{\xi}_{u}-{\xi}_{u}\in[0,\log(1+\varepsilon))\}}
\,\dd u\biggr)\,\dd e^{\alpha\overline{\xi}_{s}}\\
&&\qquad\leq e^{\alpha\overline{\xi}_{\tau(r)}}\bigl(L_{\tau(r)}+\delta
_{1}\bigr)\\
&&\qquad\quad{}-\int_{(0,\tau(r)]}L_{s}\,\dd e^{\alpha\overline{\xi
}_{s}}+\delta_{1}(e^{\alpha\overline{\xi}_{T}}-1)\\
&&\qquad\leq\int_{(0,\tau(r)]}e^{\alpha\overline{\xi}_{u}}
\,\dd L_{u}+\delta
\end{eqnarray*}
for every $r\leq t$. It follows that on $\mathcal{A}$,
\[
\frac{1}{\widehat{\mathcal{V}}(\log(1+\varepsilon))}\int
^r_{0}\mathbf{1}_{\{{M_{s}}/{X_{s}}\in[1,1+\varepsilon
)\}}\,\dd s-L^{\Theta}_{r}\leq\delta\qquad
\mbox{for
every } r\leq t,
\]
which finishes the proof of our claim.
\end{pf*}
\begin{pf*}{Proof of Theorem~\ref{exitsystem}}
Thanks to self-similarity we can suppose without loss of generality
that $x=1$. Recall that $\epsilon$ denotes the typical excursion of
the Markov process $\overline{\xi}-\xi$. For $s>0$, let $d_{s}=\inf
\{t>s\dvtx\xi_{t}=\overline{\xi}_{t}\}$. First of all, observe that the
left extrema $G$ (resp., right extrema, $D_G$) of the of the excursion
intervals complementary to the homogeneous random set $\Theta$, are
related to left extrema $g$ (resp., right extrema, $d_g$) of the
excursion intervals complementary to the regenerative set
$\{t\geq0\dvtx
\overline{\xi}_{t}-\xi_{t}=0\}$, by the relation $G=A_{g}$ for some
$g$ (resp., $D_G=A_{d_g}$ for the corresponding $d_g$). Using
this fact
and Lamperti's transformation we obtain
\begin{eqnarray*}
&&\e_{1}\biggl(\sum_{G\in\mathsf{G}}V_{G}F
\biggl(M_{G},M_{D_{G}},\biggl(\frac{M}{X}\circ\theta_{G}\biggr)\circ
k_{D_{G}}\biggr)\biggr)\\
&&\qquad=\er\biggl(\sum_{g}V_{A_{g}}F\bigl(e^{\overline{\xi}_{g}},
e^{\overline{\xi}_{g}}e^{-(\overline{\xi}_{g}-\xi_{d_{g}})},\\
&&\hspace*{92.5pt}\bigl(\exp\bigl\{(\overline{\xi}-\xi)_{\tau
(A_{g}+u)}\bigr\}, 0\leq
u\leq A_{d_{g}}-A_{g}\bigr)\bigr)\biggr)\\
&&\qquad=\er\biggl(\sum_{g}V_{A_{g}}F\bigl(e^{\overline{\xi}_{g}},\bigl(
e^{\overline{\xi}_{g}}e^{-\epsilon({d_g})}, \exp\bigl\{\epsilon
\bigl({\tau_{\epsilon}(u/e^{\alpha\overline{\xi}_{g}}
)}\bigr)\bigr\},\\
&&\hspace*{180.9pt} 0\leq u\leq e^{\alpha\overline{\xi}_{g}}A^{\{
\epsilon\}}_{d_g}\bigr)\bigr)\circ\theta_{g}\biggr),
\end{eqnarray*}
where we denote by $A^{\{\epsilon\}}_{t}=\int^t_{0}e^{\alpha\epsilon
(u)}\,\dd u$, $t\geq0$, and $\tau_{\epsilon}$ the inverse of $A^{\{
\epsilon\}}$. By the compensation formula from the excursion theory of
Markov processes
we deduce that the right-hand side in the latter equation equals
\begin{eqnarray*}
\er\biggl(\int^\infty_{0}\dd L_{s}V_{A_{s}}\overline
{n}\bigl(F\bigl(e^{\overline{\xi}_{s}}, e^{\overline{\xi
}_{s}}e^{-\epsilon({\zeta})}, \bigl(\exp\bigl\{\epsilon
\bigl({\tau_{\epsilon}(u/e^{\alpha\overline{\xi}_{s}}
)}\bigr)\bigr\}, 0\leq u\leq e^{\alpha\overline{\xi}_{s}}A^{\{
\epsilon\}}_{\zeta}\bigr)\bigr)\bigr)\biggr).
\end{eqnarray*}
Finally, using again Lamperti's transformation and the fact that
$L_{\tau}$ has as support the homogeneous random set $\Theta$, we get
that the previous expectation is equal to
\begin{eqnarray*}
&&\er\biggl(\int^\infty_{0}\dd L_{\tau(s)}V_{s}\overline
{n}\bigl(F\bigl(e^{\overline{\xi}_{\tau(s)}}, e^{\overline{\xi
}_{\tau(s)}}e^{-\epsilon({\zeta})}, \bigl(\exp\bigl\{\epsilon
\bigl({\tau_{\epsilon}\bigl(u/e^{\alpha\overline{\xi}_{\tau
(s)}}\bigr)}\bigr)\bigr\},\\
&&\hspace*{197.3pt} 0\leq u\leq e^{\alpha\overline{\xi
}_{\tau(s)}}A^{\{\epsilon\}}_{\zeta}\bigr)\bigr)\bigr)
\biggr)\\
&&\qquad=\e_{1}\biggl(\int^\infty_{0}\dd L^{\Theta
}_{s}V_{s}N^{M_{s}}(F)\biggr)\\
&&\qquad=\e_{1}\biggl(\int^\infty_{0}\dd L^{\Theta
}_{s}V_{s}N^{X_{s}}(F)\biggr),
\end{eqnarray*}
which finishes the proof.
\end{pf*}

%
\section{\texorpdfstring{Proof of Theorem \protect\ref{teo:upwardladderp}}{Proof of Theorem 2}}
\label{section:Pthm2}
The following lemma will be helpful throughout the sequel. Recall that
the process $A$ was defined in (\ref{A}).
\begin{lemma}\label{keylemma}The following equality holds $\pr$-a.s.:
%
%
\begin{equation}
A_{L^{-1}_{t}}=\int_{(0,t]}\exp\{\alpha{\mathrm h}_{s-}\}\,\dd
Y_{s},\qquad t\geq0,
\end{equation}
where we recall that $Y_{t}:=at+\sum_{u\leq t}\int
^{L^{-1}_{u}-L^{-1}_{u-}}_{0}\exp\{\alpha(\xi_{s+L^{-1}_{u-}}-\xi
_{L^{-1}_{u-}})\}\,\dd s$, $t\geq0$. The process $(L^{-1},h,Y)$
is a L\'{e}vy process with increasing coordinates and, when $\alpha
>0$, $\er(Y_{1})<\infty$.
\end{lemma}
\begin{pf}
Indeed, it follows by decomposing the interval $[0, L^{-1}_{t}]$ into
the excursion intervals that
\begin{eqnarray*}
A_{L^{-1}_{t}}&=&\int^{L^{-1}_{t}}_{0}\exp\{\alpha\xi_{s}\}\,\dd s\\[-2pt]
&=&\int^{L^{-1}_{t}}_{0}\exp\{\alpha\xi_{s}\}\mathbf{1}_{\{
\overline{\xi}_{s}-\xi_{s}=0\}}\,\dd s+\sum_{u\leq t}\int
^{L^{-1}_{u}}_{L^{-1}_{u-}}\exp\{\alpha\xi_{s}\}\,\dd s\\[-2pt]
&=&a\int^{L^{-1}_{t}}_{0}\exp\{\alpha\xi_{s}\}\,\dd L_{s}+\sum
_{u\leq t}e^{\alpha\xi_{L^{-1}_{u-}}}\int
^{L^{-1}_{u}-L^{-1}_{u-}}_{0}\exp\{\alpha(\xi_{s+L^{-1}_{u-}}-\xi
_{L^{-1}_{u-}})\}\,\dd s\\[-2pt]
&=&a\int^{t}_{0}\exp\{\alpha\xi_{L^{-1}_{u}}\}\,\dd u+\sum
_{u\leq t}e^{\alpha\xi_{L^{-1}_{u-}}}\int
^{L^{-1}_{u}-L^{-1}_{u-}}_{0}\exp\{\alpha(\xi_{s+L^{-1}_{u-}}-\xi
_{L^{-1}_{u-}})\}\,\dd s\\[-2pt]
&\geq& a\int^{t}_{0}\exp\{\alpha\xi_{L^{-1}_{u}}\}\,\dd u\\[-2pt]
&&{}+e^{(\alpha\vee
0)\xi_{L^{-1}_{t-}}}\sum_{u\leq t}\int
^{L^{-1}_{u}-L^{-1}_{u-}}_{0}\exp\{\alpha(\xi_{s+L^{-1}_{u-}}-\xi
_{L^{-1}_{u-}})\}\,\dd s.
\end{eqnarray*}
On account of the fact that for $t>0$, $A_{L^{-1}_{t}}<\infty$ a.s. we
infer that
\[
\sum_{u\leq t}\int^{L^{-1}_{u}-L^{-1}_{u-}}_{0}\exp\{\alpha(\xi
_{s+L^{-1}_{u-}}-\xi_{L^{-1}_{u-}})\}\,\dd s<\infty,\qquad
\mbox{a.s.}
\]
It follows that the subordinator $Y$ is well defined. The former
calculations and the fact that $\xi_{L^{-1}}$ has countably many
discontinuities imply that
\begin{eqnarray*}
A_{L^{-1}_{t}}&=&a\int^{t}_{0}\exp\{\alpha\xi_{L^{-1}_{u-}}\}\,\dd u+\sum
_{u\leq t}e^{\alpha\xi_{L^{-1}_{u-}}}\int
^{L^{-1}_{u}-L^{-1}_{u-}}_{0}\exp\{\alpha(\xi_{s+L^{-1}_{u-}}-\xi
_{L^{-1}_{u-}})\}\,\dd s\\[-2pt]
&=&\int_{(0,t]}\exp{\alpha\xi_{L^{-1}_{s-}}}\,\dd Y_{s}.
\end{eqnarray*}
The fact that $(L^{-1}, h,Y)$ has independent and stationary increments
is a~consequence of the fact that these processes can be explicitly
constructed in terms of the Poisson point process of excursions of $\xi
$ from its past supremum. To prove that when $\alpha>0$, $Y$ has
finite mean we use the compensation formula for Poisson point processes
\begin{eqnarray*}
\er(Y_{t})&=&at+\er\biggl(\sum_{u\leq t}\int
^{L^{-1}_{u}-L^{-1}_{u-}}_{0}\exp\{\alpha(\xi_{s+L^{-1}_{u-}}-\xi
_{L^{-1}_{u-}})\}\,\dd s\biggr)\\[-2pt]
&=&at+\er\biggl(\int^t_{0} \overline{n}\biggl(\int^{\zeta}_{0}\exp
\{-\alpha\epsilon(s)\}\,\dd s\biggr)\, \dd u\biggr)\\[-2pt]
&=&t\biggl(\int^\infty_{0}e^{-\alpha y} \widehat{\mathcal
{V}}(\dd y)\biggr)=\frac{t}{\phi_{\hat{{\mathrm h}}}(\alpha
)}<\infty,
\end{eqnarray*}
where the third and fourth identities follow from the definition of
$\widehat{\mathcal{V}}$ in (\ref{dualrenewal}) and the fact that
this is the potential measure of the downward ladder height
subordinator~$\widehat{{\mathrm h}}$, whose Laplace exponent is given by
$\phi_{\hat{{\mathrm h}}}$.
\end{pf}
\begin{pf*}{Proof of Theorem~\ref{teo:upwardladderp}} On account of
self-similarity we may prove the result under the assumption that $X_0
= 1$. We denote by $\tau_{\mathrm h}$ the time change induced by
Lamperti's transformation when applied to the ladder height
subordinator ${\mathrm h}$,
\[
\tau_{\mathrm h}(t)=\inf\biggl\{s>0\dvtx\int^{s}_{0}e^{\alpha{\mathrm
h}_{u}}\,\dd u>t\biggr\},\qquad t\geq0.
\]
By making a change of variables in the definition of the process
$L^{\Theta}$ we observe the relation
\[
L^{\Theta}_{t}=
\int^{L_{\tau(t
)}}_{0}\exp\{\alpha{\mathrm h}_{s}\}\,\dd s,\qquad t\geq0.
\]
So that its right-continuous inverse is given by $R_{t}=
A_{L^{-1}_{\tau_{\mathrm h}(t
)}}$, $t\geq0$. From Lem\-ma~\ref{keylemma}, it has the property that
its increments are given by
\begin{eqnarray*}
R_{t+s}-R_{t}&=&
\int_{(\tau_{\mathrm h}(t
),\tau_{\mathrm h}((t+s)
)]}\exp\{\alpha{\mathrm h}_{u-}\}\,\dd Y_{u}\\
&=&z^{\alpha}\int_{(0,\widetilde{\tau}_{\mathrm h}(s/z^{\alpha})]}\exp
\{\alpha\widetilde{{\mathrm h}}_{u-}\}\,\dd\widetilde
{Y}_{u}\qquad z=
\exp\bigl\{{\mathrm h}_{\tau_{\mathrm h}(t
)}\bigr\}\\
&=&H^{\alpha}_{t}\widetilde{R}_{s/H^{\alpha}_{t}},
\end{eqnarray*}
where\vspace*{2pt} $\widetilde{{\mathrm h}}_{u}={\mathrm h}_{u+\tau_{\mathrm h}(t
)}-{\mathrm h}_{\tau_{\mathrm h}(t
)}, \widetilde{Y}_{u}=Y_{u+\tau_{\mathrm h}(t
)}-Y_{\tau_{\mathrm h}(t
)}$, and given that $\tau_{{\mathrm h}}(t)$ is a~stopping time in the
filtration $\mathcal{H}_{v}:=\sigma\{({L^{-1}_{u}, \mathrm h}_{u},
Y_{u}), 0\leq u\leq v\}$, $v\geq0$, it follows that the process
$\widetilde{R}$ is a copy of $R$ issued from $0$ and independent
of~$\mathcal{H}_{\tau_{{\mathrm h}}(t)}$.

On\vspace*{1pt} account of the fact that $H$ is obtained by time changing $X$, which
is a strong Markov process, by the right-continuous inverse of the
additive functional $L^{\Theta}$, it follows by standard arguments
that $H$ is a strong Markov process. Hence the couple $(R,H)$ is a
Markov process and
the process $H$ is $1/\alpha$-self-similar.

Define $B_{H,t}=\int^t_{0}H^{-\alpha}_{s}\,\dd s$. Then we have
the following equalities:
\begin{eqnarray*}
B_{H,t}=\int_{(0,R_{t}]}X^{-\alpha}_{s}\,\dd L^{\Theta}_{s}=
\int_{(0,R_{t}]}\dd\widetilde{L}_{s}=\widetilde{L}_{R_{t}},\qquad
t\geq0,
\end{eqnarray*}
where we recall that $\widetilde{L}_t=L_{\tau(t)}$.
Denote by $C_{H,t}, t\geq0$, the inverse of the functional $B_{H,\cdot
}$. We obtain from our previous calculations that
\[
B_{H,t}=L_{\tau(A_{L^{-1}_{\tau_{{\mathrm h}}(t
)}})}=\tau_{{\mathrm h}}(t
),\qquad C_{H,t}=
\int^t_{0}e^{\alpha{\mathrm h}_{s}}\,\dd s,\qquad R_{C_{H,t}}=
A_{L^{-1}_{t}},\qquad t\geq0.
\]
These relations allow us to ensure that the process obtained by time
changing the process $\log(H_{\cdot}/H_{0})$ by $C_{H}$ is the
process ${\mathrm h}_{t}=\xi_{L^{-1}_{t}}, t\geq0$. We recall that
because of Lamperti's transformation
\[
\tau(t
)=\int^t_{0}X^{-\alpha}_{s}\,\dd s,\qquad K_{t}=\tau(
R_{t})=L^{-1}_{\tau_{h}(t
)},\qquad t\geq0.
\]
Hence the representation obtained in Lemma~\ref{keylemma} allows us to
ensure that $(K, R, H)$ time changed by the inverse of $C_{H,\cdot}$
equals $\{(L^{-1}_{t},
\int_{(0,t]}e^{\alpha{\mathrm h}_{s-}}\,\dd Y_{s},
e^{{\mathrm h}_{t}}),\break t\geq0\}$.

To prove that the process $(R,H)$ is a Feller process in $[0,\infty
)\times(0,\infty)$, we observe first that for $t,s\geq0$
\begin{eqnarray*}
e^{{\mathrm h}_{\tau_{{\mathrm h}}(t
)}
}\bigl(Y_{\tau_{{\mathrm h}}(t+s
)}-Y_{\tau_{{\mathrm h}}(t
)}\bigr)&\leq& R_{t+s}-R_{t}\leq
e^{{\mathrm h}_{\tau_{{\mathrm h}}(t+s
)}}\bigl(Y_{\tau_{{\mathrm h}}(t+s
)}-Y_{\tau_{{\mathrm h}}(t
)}\bigr), \\
H_{t+s}-H_{t}&=&
e^{{\mathrm h}_{\tau_{{\mathrm h}}(t
)}}\bigl(e^{{\mathrm h}_{\tau_{{\mathrm h}}(t+s
)}-{\mathrm h}_{\tau_{{\mathrm h}}(t
)}}-1\bigr),
\end{eqnarray*}
where $X_{0}=x=H_{0}$. By construction $\tau_{{\mathrm h}}$ is a
continuous functional, and thus if $(t_{n})_{n\geq0}$ is a convergent
sequence of positive reals with limit $t$, then the sequence of
stopping times $(\tau_{{\mathrm h}}(t_{n}))_{n\geq0}$ converges
a.s. to $\tau_{{\mathrm h}}(t)$. So, the latter inequalities together with
the right-continuity and quasi-left continuity of the L\'{e}vy process
$({\mathrm h},Y)$ (see, e.g., Proposition I.7 in~\cite{Be}) imply that if
$t_{n}\uparrow t$ or $t_{n}\downarrow t$, then
\[
(R_{t_{n}},H_{t_{n}})\xrightarrow{n\to\infty}
(R_{t},H_{t}), \qquad\p_{1}\mbox{-a.s.}
\]
Hence, the Feller
property is a simple consequence of these facts and the scaling
property.
\end{pf*}
\begin{remark}
For further generality we can chose to time change $X$ with the
right-continuous inverse of the additive functional
\[
L^{\Theta,\beta}_{t}:=\int^t_{0}X^{\beta}_{s}\,\dd\widetilde
{L}_{s},\qquad t\geq0,
\]
for some $\beta\in\re\setminus\{0\}$, fixed. In that case, we have
that the process $(K^{(\beta)},R^{(\beta)}$, $H^{(\beta)})$, defined by
\begin{eqnarray*}
R^{(\beta)}_{t}&:=&\inf\{s\geq0\dvtx L^{\Theta,\beta}_{s}>t\},\qquad
K^{(\beta)}_{t}:=\int^{R^{(\beta)}_{t}}_{0}X^{-\alpha}_{s}\,\dd s,\\
H^{(\beta)}_{t}&:=&X_{R^{(\beta)}_{t}},\qquad t\geq0,
\end{eqnarray*}
has the same law under $\p_{x}$ as the process
\[
\biggl\{\biggl(L^{-1}_{t}, x^{\alpha}\int_{(0,t]}e^{\alpha{\mathrm h}_{s-}}\,\dd
Y_{s}, xe^{{\mathrm h}_{t}}\biggr), t\geq0\biggr\},
\]
time changed by the inverse of the additive functional $(x^{\beta}\int
^t_{0}e^{\beta{\mathrm h}_{s}}\,\dd s, t\geq0)$, under~$\pr$. In
this case, the process $(R^{(\beta)}, H^{(\beta)})$ has the following
scaling property. For every $c>0$,
\[
\bigl(\bigl(c^{\alpha}R^{(\beta)}_{tc^{-\beta}},cH^{(\beta
)}_{tc^{-\beta}}\bigr), t\geq0\bigr)
\]
issued from $(x_{1},x_{2})\in\re^2_{+}$ has the same law as
$(R^{(\beta)},H^{(\beta)})$ issued from $(c^{\alpha}x_{1}$, $cx_{2})$.
The proof of these facts follows along the same lines of the proof of
Theorem~\ref{teo:upwardladderp}.
\end{remark}

\subsection{The resolvent of $X$}\label{section:res}
The main purpose of this subsection is to establish a formula for the
resolvent of $X$ in terms of the ladder process $(R,H)$. This result
will be very helpful in the proof of Theorem~\ref{thm:mainconvresult}.
\begin{corollary}\label{th:exitformulares} Let $\kappa_q(x,\mathrm
{d}y)$ be the kernel defined by
\begin{eqnarray}
\kappa_q(z,\dd y)=az^{\alpha}\delta_{z}(\mathrm
{d}y)+\overline{n}\biggl(\int^{\zeta}_{0}z^{\alpha}e^{-\alpha
\epsilon(s)}e^{-qz^{\alpha}\int^{s}_{0}e^{-\alpha\epsilon
(u)}\,\dd u}\mathbf{1}_{\{ze^{-\epsilon(s)}\in\mathrm
{d}y\}}\,\dd s\biggr),\nonumber\\
&&\eqntext{y>0, z>0.}
\end{eqnarray}

Then the $q$-resolvent of $X$, $V_{q}$, satisfies
\[
V_{q}f(x)=\e_{x}\biggl(\int^{
L^{\Theta}_{T_0}
}_{0} H^{-\alpha}_{t}e^{-qR_{t}}\kappa_q(H_{t},f)\,\dd t\biggr)
\]
for any measurable function $f\dvtx\mathbb{R}_+\to\mathbb{R}_+$ and $x>0$.
\end{corollary}
\begin{pf}
The proof is a consequence of Proposition~\ref{prop:localtime},
Theorem~\ref{exitsystem}, the identity
%
%
\begin{eqnarray}\quad
V_{q}f(x):\!&=&\e_{x}\biggl(\int^\infty_{0}e^{-qt}f(X_{t})\,\dd t\biggr)\nonumber\\
&=&\e_{x}\biggl(\int^\infty_{0}e^{-qs}\mathbf{1}_{\{X_{s}=M_{s}\}
}f(X_{s})\,\dd s\biggr)\\
&&{} +\e_{x}\biggl(\sum_{G\in\mathsf{G}}e^{-q G}\biggl(\int
^{D-G}_{0} e^{-qs}f\biggl(M_{G}\biggl(\frac{M_{G}}{X_{s+G}}
\biggr)^{-1}\biggr)\biggr)\,
\dd s\biggr)\nonumber
\end{eqnarray}
(with the notation of Theorem~\ref{exitsystem}) and that $L^{\Theta
}_{T_{0}}=\int^{L_{\infty}}_{0}e^{\alpha{\mathrm h}_{s}}\,\dd s$
under $\pr$, since~$L_{\infty}$ is the lifetime of ${\mathrm h}$.
\end{pf}

Observe that for $q=0$, the kernel $\kappa_{0}$ becomes
\begin{eqnarray*}
\kappa_{0}(z,f)&=&az^{\alpha}f(z)+\overline{n}\biggl(\int^{\zeta
}_{0}z^{\alpha}e^{-\alpha\epsilon(s)}f\bigl(ze^{-\epsilon
(s)}\bigr)\,\dd s\biggr)\\
&=&\int_{[0,\infty[}(ze^{-y})^{\alpha
}f(ze^{-y})\widehat{\mathcal{V}}(\dd y),
\end{eqnarray*}
where the second identity is a consequence of (\ref{dualrenewal}).
Recall that $H$ time changed by the inverse of $\int^t_{0}H^{-\alpha
}_{s}\,\dd s$, $t\geq0$, is equal to $e^{{\mathrm h}}$. As a
consequence we obtain that the $0$-resolvent for $X$ is given by the formula
\[
V_{0}f(x)=\iint_{[0,\infty)\times[0,\infty)}x^{\alpha}e^{\alpha
(z-y)}f(xe^{z-y})\mathcal{V}(\dd z)\widehat{\mathcal
{V}}(\dd y),
\]
where $\mathcal{V}(\dd y)$ denotes the renewal measure for
$\mathrm{h}$, that is to say
\[
\mathcal{V}(\dd y)=\er\biggl(\int^{L_\infty}_{0}\mathbf
{1}_{\{\mathrm{h}(t)\in\dd y\}}\,\dd t\biggr),\qquad
y\geq0.
\]

\section{Applications to entrance laws and weak convergence}
\label{sect:appli}

We will assume hereafter that $\alpha>0$ and $X$ is such that $\limsup
_{t\to\infty}X_{t}=\infty$ a.s. which is equivalent to assume that,
for the underlying L\'{e}vy process $\xi$, it holds that $\limsup
_{t\to\infty}\xi_{t}=\infty$, a.s.
\begin{lemma}\label{keylemma:bis}
For every $t>0$, the following equality in law holds:
%
%
\begin{eqnarray}\label{eq:dualityrel}
&&\bigl(\bigl({\mathrm h}_{t}-{\mathrm h}_{(t-s)-}, 0\leq s\leq t\bigr),
e^{-\alpha{\mathrm h}_{t}}A_{L^{-1}_{t}}\bigr)\nonumber\hspace*{-35pt}\\[-8pt]\\[-8pt]
&&\qquad \stackrel{(d)}{=}\biggl(({\mathrm h}_{s},
0\leq s\leq t), \int_{(0,t]}e^{-\alpha{\mathrm h}_{r-}}\,\dd Y_{r}+\sum
_{u\leq t}(e^{-\alpha\Delta{\mathrm h}_{u}}-1)e^{\alpha
{\mathrm h}_{u-}}\Delta Y_{u}\biggr).
\nonumber\hspace*{-35pt}
\end{eqnarray}
Furthermore, if $\alpha>0$ and $\er({\mathrm h}_{1})<\infty$, then the
stochastic process
\[
W_{t}:=\int_{(0,t]}e^{-\alpha{\mathrm h}_{r-}}\,\dd Y_{r}+\sum
_{u\leq t}(e^{-\alpha\Delta{\mathrm h}_{u}}-1)e^{- \alpha{\mathrm
h}_{u-}}\Delta Y_{u}
\]
converges a.s., as $t\to\infty$, to a random variable
$\widetilde{I}$ satisfying
\[
\widetilde{I}:=\int_0^{\infty} \exp\{-\alpha{\mathrm h}_s\}\,\dd
Y_{s-}\stackrel{(d)}{=}\int_0^{\infty} \exp\{-\alpha\xi^{\uparrow
}_s\}\,\dd s.
\]
\end{lemma}
\begin{pf}
We start by proving the time reversal property described in~(\ref
{eq:dualityrel}). On the one hand the duality lemma for L\'{e}vy
processes implies that for $t>0$
\[
\bigl(\bigl({\mathrm h}_{t}-{\mathrm h}_{(t-s)-}, Y_{t}-Y_{(t-s)-}\bigr), 0\leq s\leq
t\bigr)\stackrel{(\dd)}{=}\bigl(({\mathrm h}_{s},Y_{s}), 0\leq
s\leq t\bigr).
\]
It follows, making a change of variables of the form $t-u$ together
with the above identity in law, that
\begin{eqnarray*}
e^{-\alpha{\mathrm h}_{t}}\int_{(0,t]}e^{\alpha{\mathrm h}_{s-}}\,\dd
Y_{s}&=&\int_{(0,t]}e^{-\alpha({\mathrm h}_{t}-{\mathrm h}_{(t-u)-})}\,
\dd
(Y_{t}-Y_{t-u})\\
&\stackrel{(\mathrm{d})}{=}&\int_{(0,t]}e^{-\alpha{\mathrm h}_{s}}\,\dd Y_{s-}\\
&=&\int_{(0,t]}e^{-\alpha{\mathrm h}_{r-}}\,\dd Y_{r}+\sum_{u\leq
t}(e^{-\alpha\Delta{\mathrm h}_{u}}-1)e^{-\alpha{\mathrm h}_{u-}}\Delta Y_{u}.
\end{eqnarray*}
A similar identity for general L\'{e}vy processes can be found in
\cite{LindnerMaller2005}. Observe that the process $T$ defined by
$T_{t}=Y_{t}+\sum_{s\leq t}(e^{-\alpha\Delta{\mathrm h}_{s}}-1)\Delta
Y_{s}$, $t\geq0$, is a~L\'{e}vy process and for any $t>0$,
\[
\int_{(0,t]}e^{-\alpha{\mathrm h}_{r-}}\,\dd Y_{r}+\sum_{u\leq
t}(e^{-\alpha\Delta{\mathrm h}_{u}}-1)e^{-\alpha{\mathrm h}_{u-}}\Delta
Y_{u}=\int_{(0,t]}e^{-\alpha{\mathrm h}_{s-}}dT_{s}.
\]
To finish the proof of the convergence of $W_t$, we should verify that,
when $\alpha>0$ and $\er({\mathrm h}_{1})<\infty$, the right-hand side
of the above equality converges a.s. Thanks to the assumption that $\er
({\mathrm h}_{1})<\infty$, according to Theorem 2.1 and Remark 2.2 in
\cite{LindnerMaller2005} it suffices to verify that $\er(\log
^+(T_{1}))<\infty$. Indeed, given that $T_{1}\leq Y_{1}$ a.s. it
follows from Lemma~\ref{keylemma} that under our assumptions $\er
(Y_{1})<\infty$, and thus $\er(\log^+(Y_{1}))<\infty$ which in turn
implies that $\er(\log^+(T_
{1}))<\infty$.

Finally the identity in law between $\widetilde{I}$ and $\int^\infty
_{0}\exp\{-\alpha\xi^{\uparrow}_{s}\}\,\dd s$ follows from the
Doney--Tanaka path construction of the process conditioned to stay
positive. Roughly speaking, the latter allows us to construct the
process $\xi^{\uparrow}$ by pasting, at the level of the last
supremum, the excursions of $\xi$ from its past supremum, reflected
and time reversed. More precisely, for $t\geq0$, let
\begin{eqnarray*}
g_{t}&=&\sup\{s<t\dvtx\overline{\xi}_{s}-\xi_{s}=0\},\qquad d_{t}=\inf\{
s>t\dvtx\overline{\xi}_{s}-\xi_{s}=0\},
\\
\widetilde{\mathcal{R}}_{t}&=&\cases{
(\overline{\xi}-\xi)_{(d_{t}+g_{t}-t)-}, &\quad if $d_{t}-g_{t}>0$,\cr
0, &\quad if $d_{t}-g_{t}=0$,}\qquad
\mathcal{R}_{t}=\overline{\xi}_{d_{t}}+\widetilde{\mathcal{R}}_{t}.
\end{eqnarray*}
The process $\{\mathcal{R}_{t}, t\geq0\}$, under $\pr$ has the same
law as $(\xi^{\uparrow},\pr^{\uparrow})$. See~\cite{doney}, Section~8.5.1,
for a proof of this result. Taking account of the Doney--Tanaka
construction, we may proceed as in the proof of Lemma~\ref{keylemma}
and obtain that
\begin{eqnarray*}
&&\int^\infty_{0}\exp\{-\alpha\xi^\uparrow_{s}\}\,\dd s\\
&&\qquad\stackrel{(\dd)}{=}\int^{\infty}_{0}\exp\{-\alpha\overline
{\xi}_{d_{s}}\}\mathbf{1}_{\{d_{s}=g_{s}\}}\,\dd s\\
&&\qquad\quad{}
+\sum_{t\in \widetilde{\mathsf{G}}}\int^{d_{t}}_{g_{t}}\exp\bigl\{-\alpha\bigl(\overline
{\xi}_{d_{s}}+(\overline{\xi}-\xi)_{(d_{s}+g_{s}-s)-}\bigr)\bigr\}\,\dd s,
\end{eqnarray*}
where $\widetilde{\mathsf{G}}$ denotes the left extrema of the
excursion intervals of $\xi$ from $\overline{\xi}$.
We recall that $g_{s}, d_{s}$, remain constant along the excursion
intervals, that $d_{s}=L^{-1}_{L_{s}}$, $g_{s}=L^{-1}_{L_{s}-}$, and
that $d_{s}=g_{s}$ if and only if $s$ belongs to the interior of the
random set of points where $\overline{\xi}-\xi$, takes the value
$0$. Using these facts we obtain the following identities:
\begin{eqnarray*}
&&\int^\infty_{0}\exp\{-\alpha\xi^{\uparrow}_{s}\}\,\dd s\\[-1pt]
&&\qquad\stackrel{(\dd)}{=}\int^\infty_{0}\exp\{-\alpha\overline
{\xi}_{s}\}\mathbf{1}_{\{\xi_{s}=\overline{\xi}_{s}\}}\,\dd s\\[-1pt]
&&\qquad\quad{} + \sum
_{t\geq0}\int^{L^{-1}_{t}}_{L^{-1}_{t-}}\exp\bigl\{
-\alpha\bigl(\overline{\xi}_{L^{-1}_{t}}+(\overline{\xi
}-\xi)_{{L^{-1}_{t}+L^{-1}_{t-}-s}}\bigr)\bigr\}\,\dd s\\[-1pt]
&&\qquad=a\int^\infty_{0}\exp\{-\alpha\overline{\xi}_{s}\}\,\dd L_{s}\\[-1pt]
&&\qquad\quad{}+\sum
_{t\geq0}e^{-\alpha\overline{\xi}_{L^{-1}_{t}}}\int
^{L^{-1}_{t}}_{L^{-1}_{t-}}\exp\{-\alpha(\overline{\xi
}_{L^{-1}_{t-}}-\xi_{L^{-1}_{t}+L^{-1}_{t-}-s})\}\,\dd s\\[-1pt]
&&\qquad=a\int^\infty_{0}\exp\{-\alpha\overline{\xi}_{s}\}\,\dd L_{s}\\[-1pt]
&&\qquad\quad{} + \sum
_{t\geq0}e^{-\alpha\overline{\xi}_{L^{-1}_{t}}}\int
^{L^{-1}_{t}-L^{-1}_{t-}}_{0}\exp\{-\alpha(\overline{\xi
}_{L^{-1}_{t-}}-\xi_{L^{-1}_{t-}+u})\}\,\dd u\\[-1pt]
&&\qquad=\int^\infty_{0}\exp\{-\alpha\overline{\xi}_{L^{-1}_{t}}\}\,
\dd Y_{t}.
\end{eqnarray*}
To conclude we should justify that
\[
\int^\infty_{0}\exp\{-\alpha\overline{\xi}_{L^{-1}_{t}}\}\,\dd
Y_{t}=\int^\infty_{0}\exp\{-\alpha\overline{\xi}_{L^{-1}_{t}}\}\,
\dd Y_{t-}, \qquad\pr\mbox{-a.s.}
\]
Indeed, we have that for every $s>0$
\[
\int_{(0,s]}\exp\{-\alpha\overline{\xi}_{L^{-1}_{t}}\}\,\dd Y_{t}=\int
_{(0,s)}\exp\{-\alpha\overline{\xi}_{L^{-1}_{t}}\}\,
\dd Y_{t-}+(Y_{s}-Y_{s-})e^{-\alpha\overline{\xi}_{L^{-1}_{s}}}.
\]
Furthermore, the last term on the right-hand side above tends to $0$ as
$s$ growths to infinity, $\pr$-a.s. because $Y_{s}$ and $\overline
{\xi}_{L^{-1}_{s}}$ have linear growth, which in turn is thanks to the
hypothesis $\er(\mathrm{h}_{1})<\infty$, and the fact $\er
(Y_{1})<\infty$ (which was established in Lemma~\ref{keylemma}). The
identity follows.
\end{pf}
\begin{pf*}{Proof of Theorem~\ref{limitladder}}
Our objective will be achieved in three main steps. First, we will
prove that the resolvent of $(R,H)$ under $\p_{x}$ has a
nondegenerate limit as $x\to0$. Second, we will deduce the
finite-dimensional convergence,\vadjust{\goodbreak} using an argument that follows along
the lines
of Bertoin and Yor's proof of the finite-dimensional convergence of a
pssMp as the starting point tends to $0$. (Our argument is essentially
a rewording of those in the three paragraphs following the second
display on page 396 of~\cite{BeY}. We include the arguments here for
sake of completeness.) We recall that the identity in law between
$\widetilde{I}$ and $\int^\infty_{0}\exp\{-\alpha\xi^{\uparrow
}_{s}\}\,\dd s$, has been proved in Lemma~\ref{keylemma:bis}. So,
to finish we will prove the formula in Theorem~\ref{limitladder} for
the limit law of $R,H$ under $\p_{x}$ as $x\to0$.

\textit{Step} 1. Observe that by construction the process
$(R,H)$ issued from $(R_{0}=y$, $H_{0}=x)$ has the same law as $(R+y,H)$
issued from $(R_{0}=0,H_{0}=x)$. We denote by $V^{R,H}_{q}, q\geq0$,
the $q$-resolvent of the process $(R,H)$, namely for continuous and
bounded $f\dvtx\re_+^{2}\to\mathbb{R}_+$,
%
%
\begin{equation}\label{resRH}
V^{R,H}_{q}f(y,x):=\e_{x}\biggl(\int^\infty
_{0}e^{-qt}f(y+R_{t},H_{t})\,\dd t\biggr),
\end{equation}
where we recall that by construction $R_{0}=0$, $\p_{x}$-a.s. As this
operator is clearly continuous in $y$ we can assume without loss of
generality that $y=0$. We will prove that 
$V^{R,H}_{q}f(0,x)$ has a nondegenerate limit as $x\to0+$.
Recall that under $\p_{x}$ the process $((R_{t},H_{t})_{t\geq
0})$ is equal in law to the process
\[
\bigl(\bigl(x^{\alpha}A_{L^{-1}_{\tau_{\mathrm h}(t/x^{\alpha})}}, xe^{{\mathrm
h}_{\tau_{\mathrm h}(t/x^{\alpha})}}\bigr), t\geq0\bigr)\qquad \mbox
{under } \pr.
\]
We will need the following identity for $f\dvtx\re^2\to\mathbb{R}_+$ measurable:
%
%
\begin{eqnarray}\label{eq:dualityuplhp}
&&\e_{x}\biggl(\int^\infty_{0}e^{-qt}f(R_{t},H_{t})\,\dd t
\biggr)\nonumber\\[-2pt]
&&\qquad=\er\biggl(\int^\infty_{0}e^{-qx^{\alpha}\int^t_{0}e^{\alpha
\mathrm{h}_{u}}\,\dd u}x^{\alpha}e^{\alpha\mathrm
{h}_{t}}f(x^{\alpha}A_{L^{-1}_{t}}, xe^{\mathrm{h}_{t}})\,\dd
t\biggr)\nonumber\\[-2pt]
&&\qquad=\int^\infty_{0} \er\biggl(\exp\biggl\{-qx^{\alpha} e^{\alpha
\mathrm{h}_{t}}\int^t_{0}e^{-\alpha(\mathrm{h}_{t}-\mathrm
{h}_{u})}\,\dd u\biggr\}\nonumber\\[-2pt]
&&\qquad\quad\hspace*{31pt}{}\times x^{\alpha}e^{\alpha\mathrm
{h}_{t}}f(x^{\alpha}e^{\alpha\mathrm{h}_{t}}e^{-\alpha\mathrm
{h}_{t}}A_{L^{-1}_{t}}, xe^{\mathrm{h}_{t}})\,\dd
t\biggr)\nonumber\\[-9pt]\\[-9pt]
&&\qquad=\int^\infty_{0} \er\biggl(\exp\biggl\{-qx^{\alpha} e^{\alpha
\mathrm{h}_{t}}\int^t_{0}e^{-\alpha\mathrm{h}_{u-}}\,\dd u\biggr\}\nonumber\\[-2pt]
&&\qquad\quad\hspace*{31pt}{}\times x^{\alpha
}e^{\alpha\mathrm{h}_{t}}f(x^{\alpha
}e^{\alpha\mathrm{h}_{t}}W_{t}, xe^{\mathrm{h}_{t}})\,\dd t\biggr)\nonumber\\[-2pt]
&&\qquad=\er\biggl( \int^\infty_{0}\exp\biggl\{-q x^{\alpha}e^{\alpha
{\mathrm h}_{\tau_{\mathrm h}(u/x^{\alpha})}}\int^{\tau_{\mathrm
h}(u/x^{\alpha
})}_{0}e^{-\alpha{\mathrm h}_{s-}}\,\dd s\biggr\}\nonumber\\[-2pt]
&&\qquad\quad\hspace*{32pt}{}\times{f}\bigl(x^\alpha
e^{\alpha{\mathrm h}_{\tau_{\mathrm h}(u/x^{\alpha})}} W_{\tau_{\mathrm
h}(u/x^{\alpha})},xe^{{\mathrm h}_{\tau_{\mathrm h}(u/x^{\alpha})}}\bigr)\,
\dd
u\biggr),\nonumber
\end{eqnarray}
where in the first equality we made the change of variables $u=\tau
_{h}(t);$ in the second we applied Fubini's theorem;\vadjust{\goodbreak}
in the third we
used the time inversion property obtained in Lemma~\ref{keylemma:bis}
(note also that the process $W$ was defined in Lemma \ref
{keylemma:bis}); finally in the fourth we used Fubini's theorem again
and applied the change of variables $t=\tau_{\mathrm{h}}(u)$.

Next we note that a consequence of Lemma~\ref{keylemma:bis} is that
under our assumptions the pair $(W_{\tau_{\mathrm h}(s)}, \int^{\tau
_{\mathrm h}(s)}_{0}e^{-\alpha{\mathrm h}_{t-}}\,{\dd}t)$ converges a.s. to
$(\widetilde{I},I_{\mathrm h})$, as $s\to\infty$. The fact that
${\mathrm h}$ is nonarithmetic and $\er({\mathrm h}_{1})<\infty$ imply that
$H_{t}t^{-1/\alpha}$ converges weakly to a nondegenerated random
variable $Z$, by the main theorem in~\cite{beC}. Using these results
we deduce that
\begin{eqnarray*}
&&\biggl(u^{-1/\alpha}e^{{\mathrm h}_{\tau(u)}},u^{-1/\alpha}e^{{\mathrm
h}_{\tau(u)}}\int^{\tau(u)}_{0}e^{-\alpha{\mathrm h}_{s-}}\,\dd s,
\bigl(u^{-1/\alpha}e^{{\mathrm h}_{\tau(u)}}\bigr)^{\alpha}W_{\tau
_{\mathrm h}(u)}\biggr)\\
&&\qquad\xrightarrow{u\to\infty}^{(\dd)}\biggl(Z, Z\int^{\infty
}_{0}e^{-\alpha{\mathrm h}_{s}}\,\dd s, Z^{\alpha}\widetilde{I}\biggr)
\end{eqnarray*}
under $\pr$. Thus if $f\dvtx\mathbb{R}_+^2\to\mathbb{R}_+$ is a
continuous and bounded function then equation (\ref{eq:dualityuplhp})
and Fatou's lemma imply that
\[
\liminf_{x\to0+}\e_{x}\biggl(\int^\infty_{0}\mathrm
{d}te^{-qt}f(R_{t},H_{t})\biggr)
\geq\int^\infty_{0}\er(\exp\{-qt Z^{\alpha}I_{\mathrm h}\}
f(tZ^{\alpha}\widetilde{I},t^{1/\alpha}Z))\,\dd t.
\]
Furthermore, let $M=\sup_{(z,y)\in\mathbb{R}_+^2}f(y,z)$ and
$f^c(z,y)=M-f(z,y)$, and apply the latter estimate to $f^c$ to get that
\begin{eqnarray*}
&&\frac{M}{q}-\limsup_{x\to0+}\e_{x}\biggl(\int^\infty
_{0}e^{-qt}f(R_{t},H_{t})\,\dd t\biggr)\\
&&\qquad=\liminf_{x\to0+}\e_{x}\biggl(\int^\infty
_{0}e^{-qt}f^c(R_{t},H_{t})\,\dd t\biggr)\\
&&\qquad\geq\int^\infty_{0}\er(\exp\{-qt Z^{\alpha}I_{\mathrm h}\}
f^c(tZ^{\alpha}\widetilde{I},t^{1/\alpha}Z))\,\dd t\\
&&\qquad= M\int^\infty_{0}\er(\exp\{-qt Z^{\alpha}I_{\mathrm h}\}
)\,\dd t\\
&&\qquad\quad{}-\int^\infty_{0}\er(\exp\{-qt Z^{\alpha}I_{\mathrm h}\}
f(tZ^{\alpha}\widetilde{I},t^{1/\alpha}Z))\,\dd t.
\end{eqnarray*}
We have thus proved the inequalities
\begin{eqnarray*}
&&\int^\infty_{0}\er(\exp\{-qt Z^{\alpha}I_{\mathrm h}\}f(tZ^{\alpha
}\widetilde{I},t^{1/\alpha}Z))\,\dd t\\
&&\qquad\leq
\liminf_{x\to0+}\e_{x}\biggl(\int^\infty
_{0}e^{-qt}f(R_{t},H_{t})\,\dd t\biggr)\\
&&\qquad\leq\limsup_{x\to0+}\e_{x}\biggl(\int^\infty
_{0}e^{-qt}f(R_{t},H_{t})\,\dd t\biggr)\\
&&\qquad\leq\frac{M(1-\er((Z^{\alpha}I_{\mathrm h}
)^{-1}))}{q}\\
&&\qquad\quad{}+\int^\infty_{0}\er(\exp\{-qt Z^{\alpha
}I_{\mathrm h}\}f(tZ^{\alpha}\widetilde{I},t^{1/\alpha}Z))\,\dd t.
\end{eqnarray*}
Next, we need to verify that
%
%
\begin{equation}\label{eq:unidad}\er((Z^{\alpha}I_{\mathrm h})^{-1})=1.
\end{equation}
To this end, observe the following duality identity:
%
%
\begin{eqnarray}\label{eq:duality0}\qquad
&&\int^\infty_{0} e^{-\lambda t}\er\biggl(\exp\biggl\{-q x^{\alpha
}e^{\alpha{\mathrm h}_{\tau_{\mathrm h}(t/x^{\alpha})}}\int^{\tau
_{\mathrm h}(t/x^{\alpha})}_{0}e^{-\alpha{\mathrm h}_{s-}}\,\dd s\biggr\}
\biggr)\,\dd t\nonumber\\[-8pt]\\[-8pt]
&&\qquad=\int^\infty_{0} e^{-q t}\er\biggl(\exp\biggl\{-\lambda x^{\alpha
}e^{\alpha{\mathrm h}_{\tau_{\mathrm h}(t/x^{\alpha})}}\int^{\tau
_{\mathrm h}(t/x^{\alpha})}_{0}e^{-\alpha{\mathrm h}_{s-}}\,\dd s\biggr\}\biggr)\,
\dd
t,\nonumber
\end{eqnarray}
valid for $\lambda, q\geq0$. The proof of this identity uses the same
arguments as in (\ref{eq:dualityuplhp}). Taking $x\to0$ and using the
dominated convergence theorem we deduce from (\ref{eq:duality0}) that
%
%
\begin{eqnarray}
\er\biggl(\frac{1}{\lambda+q Z^{\alpha}I_{h}}\biggr)&=&\int^\infty
_{0}e^{-\lambda t}\er(\exp\{-q Z^{\alpha}I_{\mathrm
{h}}t\})\,\dd t\nonumber\\
&=&\int^\infty_{0}e^{-q t}\er(\exp\{-\lambda Z^{\alpha
}I_{\mathrm{h}}t\})\,\dd t\\
&=&\er\biggl(\frac{1}{q+\lambda Z^{\alpha}I_{h}}\biggr)\nonumber
\end{eqnarray}
for $q,\lambda>0$. We may now let $q$ tend to $0+$ and apply the
monotone convergence theorem to obtain the claimed equality.

The latter facts and a change of variables lead to
\begin{eqnarray*}
\lim_{x\to0+}\e_{x}\biggl(\int^\infty
_{0}e^{-qt}f(R_{t},H_{t})\biggr)&=&\int^\infty_{0}\er(\exp\{-qt
Z^{\alpha}I_{\mathrm h}\}f(tZ^{\alpha}\widetilde{I},t^{1/\alpha
}Z))\,\dd t\\
&=&\int^\infty_{0} e^{-qt}\er\biggl(f\biggl(\frac{t\widetilde
{I}}{I_{\mathrm h}},\frac{t^{1/\alpha}}{I^{1/\alpha}_{\mathrm h}}
\biggr)\frac{1}{Z^{\alpha}I_{\mathrm h}}\biggr)\,\dd t.
\end{eqnarray*}
We have thus proved the convergence of the resolvent of $(R,H)$ under
$\p_{x}$ as \mbox{$x\to0$}.

\textit{Step} 2. We define for $f\dvtx\re_+^{2}\to\re_{+}$ measurable
%
%
\begin{eqnarray}\label{VRH0}
V^{R,H}_{q}f(y,0)&=&\lim_{x\to0+}\e_{x}\biggl(\int^\infty
_{0}e^{-qt}f(y+R_{t},H_{t})\,\dd t\biggr)\nonumber\\[-8pt]\\[-8pt]
&=&\int^\infty_{0} e^{-qt}\er\biggl(f\biggl(y+\frac{t\widetilde
{I}}{I_{\mathrm h}},\frac{t^{1/\alpha}}{I^{1/\alpha}_{\mathrm h}}
\biggr)\frac{1}{Z^{\alpha}I_{\mathrm h}}\biggr)\,\dd t,\qquad y\geq
0.\nonumber
\end{eqnarray}
Let $C_{0}$ be the space of continuous functions on $\mathbb{R}^2_+$
with limit $0$ at infinity. Observe that for every $t>0$, the stopping
time $\tau_{h}(tx^{-\alpha})$ tends to $\infty$ (resp., to $0$) as
$x\to0$ (resp., to $\infty$) and that $A_{L^{-1}_{s}}$ tends to
$\infty$ (resp., to $0$) as $s\to\infty$ (resp., $s\to0$). Moreover,
given that
\[
A_{u}/u\xrightarrow{u\to0} 1,\qquad \tau_{h}(u)/u\xrightarrow{u\to0}1
\quad\mbox{and}\quad L^{-1}_{u}/u\xrightarrow{u\to0} a,
\]
$\pr$-almost surely, it follows that
\[
x^{\alpha}A_{L^{-1}_{\tau_{h}(x^{-\alpha})}}\xrightarrow{x\to
\infty}a,
\]
$\pr$-a.s. Using these facts and the results in Theorem \ref
{teo:upwardladderp}(i), we deduce that for every function $f\in
C_{0}$, the function $x,y\mapsto\e_{x}(f(y+R_{t},H_{t})
)$ has a limit $0$ at infinity. Now, applying the result in
Theorem~\ref{teo:upwardladderp}(ii) it follows that the operator
$V^{R,H}_{q}$ maps $C_{0}$ into $C_{0}$, and furthermore the resolvent
equation holds on~$C_{0}$. It follows from (\ref{eq:unidad}), (\ref
{VRH0}) and the Feller property in Theorem~\ref{teo:upwardladderp}
that for $f\in C_{0}$, $qV^{R,H}_{q}$ converges pointwise to $f$ as
$q\to\infty$. By the discussion on page 83 in~\cite{revuzyor} it
follows that this implies the uniform convergence of $qV^{R,H}_{q}f$
toward $f$ as $q\to\infty$ for $f\in C_{0}$. Now we invoke the
Hille--Yoshida theory to deduce that
associated to the resolvent family $qV^{R,H}_{q}$ there is a unique
strongly continuous Markovian semi-group on $C_{0}$. The
finite-dimensional convergence follows.

\textit{Step} 3.
We will next establish the formula for the entrance law for $(R,H)$. On
the one hand, by the scaling property of $(R,H)$ the weak convergence of
the one-dimensional law of $(R,H)$ as the starting point of $X$ tends
to $0$, is equivalent to the weak convergence of $(t^{-1}R_{t},
t^{-1/\alpha}H_{t})$ as $t\to\infty$, under $\p_{1}$,
as $t\to\infty$. Moreover, using the self-similarity and Feller
properties of $(R,H)$ and arguing as in Sections~5.2 to 5.4 of
\cite{CPY}, one may show that the process
\[
\mathrm{OU}_{t}:=(e^{-t}R_{e^t-1}, e^{-t/\alpha}H_{e^t-1}),\qquad
t\geq0,
\]
is an homogeneous Markov process which, under our assumptions, has
a~unique invariant measure. On the other hand, it has been proved in
\cite{r2008} that for any $c>0$, the measure defined for $f\dvtx\re
_+^2\to\mathbb{R}_+$ positive and measurable by
\[
c\er\biggl(f\biggl(\frac{s\widetilde{I}}{I_{\mathrm h}},\frac
{s^{1/\alpha}}{I^{1/\alpha}_{\mathrm h}}\biggr)\frac{1}{I_{\mathrm h}}\biggr)
\]
is an entrance law for the process $(R,H)$. It readily follows that
this measure, taking $s=1$, is an invariant measure for the process
$\mathrm{OU}$. Furthermore, the constant $c$ can be chosen so that the
measure defined above is a probability measure. This follows as a
consequence of the fact that under our assumptions, $\er
(I^{-1}_{\mathrm h})=\alpha\mu_{+}<\infty$ (see, e.g.,~\cite{beC} and
\cite{BeY}).
Therefore, we conclude that
\[
\frac{1}{\alpha\mu_{+}}\er\biggl(f\biggl(\frac{s\widetilde
{I}}{I_{\mathrm h}},\frac{s^{1/\alpha}}{I^{1/\alpha}_{\mathrm h}}
\biggr)\frac{1}{I_{\mathrm h}}\biggr)=\er\biggl(f\biggl(\frac{s\widetilde
{I}}{I_{\mathrm h}},\frac{s^{1/\alpha}}{I^{1/\alpha}_{\mathrm h}}
\biggr)\frac{1}{Z^{\alpha}I_{\mathrm h}}\biggr)
\]
for every positive and measurable function $f$. This implies in
particular that
\[
\er(Z^{-\alpha}| I_{\mathrm h}, \widetilde{I})=\frac{1}{\alpha\mu_{+}}.
\]
Which finishes the proof of Theorem~\ref{limitladder}.
\end{pf*}

\subsection{\texorpdfstring{Proof of Theorem \protect\ref{thm:mainconvresult}: Finite-dimensional convergence}
{Proof of Theorem 4: Finite-dimensional convergence}}\label{proofofthmconvergence}

We will prove this result in two main steps, first of all we will prove
the convergence in the finite-dimensional sense and then in Section
\ref{weakconvergence} we will prove that the convergence holds in the
Skorohod's sense.

As in the proof of Theorem~\ref{limitladder}, the most important tool
to prove convergence of finite-dimensional distributions is to
establish the convergence of the resolvent of~$X$ as $x$ tends to $ 0$.
One may then appeal to reasoning along the lines the proof~of Theorem 1
in~\cite{BeY} (see also the second step in the proof of our
Theorem~\ref{limitladder}). Note that, while Bertoin and Yor
\cite{BeY} require that the underlying L\'{e}vy process drifts to
$\infty$,
we may circumvent the use of this condition as we are able to write the
resolvent of $X$ in terms of the process $(R,H)$. We omit the details.
We will finish the proof by establishing the formula for the entrance law.

First of all we recall that in Corollary~\ref{th:exitformulares} we
established that the $q$-resolvent of $X$ is given by
%
%
\begin{equation}\label{eqresolvent}
V_{q}f(x)=\int^\infty_{0}\e_{x}(H^{-\alpha
}_{t}e^{-qR_{t}}\kappa_{q}(H_{t},f))\,\dd t
\end{equation}
for every $f\dvtx\mathbb{R}_+\to\mathbb{R}_+$ measurable and bounded.
It follows that for $q>0$,
%
%
\begin{equation}\label{eq:resolvent1}
\e_{x}\biggl(\int^{\infty}_{0} e^{-q R_{t}}H^{-\alpha}_{t}\kappa
_{q}(H_{t},1)\,\dd t\biggr)=\frac{1}{q},\qquad x>0.
\end{equation}
To prove the convergence of the resolvent $V_{q}$ it will be useful to
know that the mapping $x\mapsto x^{-\alpha}\kappa_{q}(x,1)$, for
$x>0$ defines a decreasing, continuous and bounded function. Indeed for
every $x>0$
\begin{eqnarray*}
x^{-\alpha}\kappa_{q}(x,1)&=&a+(qx^{\alpha})^{-1}\overline{n}
\biggl(1-\exp\biggl\{-qx^{\alpha}\int^\zeta_{0}e^{-\alpha\epsilon
(s)}\,\dd s\biggr\}\biggr)\\
&=&a+\int^\infty_{0}\,\dd y\overline{n}\biggl(\int^\zeta
_{0}e^{-\alpha\epsilon(s)}\,\dd s>y\biggr) e^{-qx^{\alpha}y}.
\end{eqnarray*}
Given that for $\alpha>0$
\[
\lim_{x\to0+}\frac{\kappa_{q}(x,1)}{x^{\alpha}}=a+\overline
{n}\biggl(\int^\zeta_{0}e^{-\alpha\epsilon(s)}\,\dd s\biggr)=\er
(Y_{1})<\infty,
\]
the claim follows. The previous limit together with the dominated
convergence theorem imply that if $f\dvtx\mathbb{R}_+\to\re$ is a
continuous and bounded function
then so is $x\mapsto x^{-\alpha}\kappa_{q}(x,f)$, for $x>0$. So,
Theorem~\ref{limitladder} implies that if $f$ is a~continuous and
bounded function, then
\[
\lim_{x\to0+}\e_{x}(e^{-qR_{t}}H^{-\alpha}_{t}\kappa
_{q}(H_{t},f))=\e^{R,H}_{0+}(e^{-qR_{t}}H^{-\alpha
}_{t}\kappa_{q}(H_{t},f)).
\]
This limit result together with Fatou's lemma imply that
%
%
\begin{equation}\label{liminfVq}
\liminf_{x\to0+}V_{q}f(x)\geq\int^\infty_{0}\e^{R,H}_{0+}
(e^{-qR_{t}}H^{-\alpha}_{t}\kappa_{q}(H_{t},f))\,\dd t.
\end{equation}
To determine the upper limit of $V_{q}f$ we first claim that
%
%
\begin{equation}\label{equalityres1}\int^\infty_{0}\e
^{R,H}_{0+}(e^{-qR_{t}}H^{-\alpha}_{t}\kappa_{q}(H_{t},1)
)\,\dd t=\frac{1}{q}.
\end{equation}
To see why the above claim is true, note that $\p_{0+}(R_{t}\in\dd s,
H_{t}\in\dd y)$ is an entrance law for the semigroup of the
self-similar Markov process $(R,H)$, and hence for any $\epsilon>0$
%
%
\begin{eqnarray}\label{star}
&&\int^\infty_{0}\e^{R,H}_{0+}(e^{-qR_{t}}H^{-\alpha
}_{t}\kappa_{q}(H_{t},1))\,\dd t\nonumber\\
&&\qquad=\int^{\varepsilon}_{0}\e^{R,H}_{0+}(e^{-qR_{t}}H^{-\alpha
}_{t}\kappa_{q}(H_{t},1))\,\dd t\nonumber\\
&&\qquad\quad{} +\int_{\mathbb{R}_+\times(0,\infty)}\p
^{R,H}_{0+}(R_{\varepsilon}\in du,H_{\varepsilon}\in
dx)\nonumber\\[-8pt]\\[-8pt]
&&\qquad\quad\hspace*{11pt}{}\times\int^\infty
_{0}\e_{x}\bigl(e^{-q(u+R_{t})}H^{-\alpha}_{t}\kappa
_{q}(H_{t},1)\bigr)\,\dd t\nonumber\\
&&\qquad \leq\er(Y_1)\int^{\varepsilon}_{0}\e^{R,H}_{0+}
(e^{-qtR_{1}})\,\dd t +\frac{1}{q}\e^{R,H}_{0+}(e^{-q
R_{\varepsilon}})\nonumber\\
&&\qquad \leq\er(Y_1)\varepsilon+\frac{1}{q}\e^{R,H}_{0+}(e^{-q
\varepsilon R_{1}}),\nonumber
\end{eqnarray}
where the first inequality comes from (\ref{eq:resolvent1}) and the
fact that $x^{-\alpha}\kappa_{q}(x,1)\leq\er(Y_{1})<\infty$,
$\forall x\geq0$.
Identity (\ref{equalityres1}) is obtained by taking $\varepsilon\to
0+$ in (\ref{star}) and combining the resulting inequality together
with an analogous lower bound which can be obtained in a similar way.

By applying the inequality (\ref{liminfVq}) to the function
$f^c(x)=\sup_{z>0}f(z)-f(x)$, $x>0$, and using the latter identity we
obtain that
%
%
\begin{equation}\label{eqlimsupVq}\limsup_{x\to0+}V_{q}f(x)\leq\int
^\infty_{0}\e^{R,H}_{0+}(e^{-qR_{t}}H^{-\alpha}_{t}\kappa
_{q}(H_{t},f))\,\dd t.
\end{equation}
We have therefore proved that for every continuous and bounded function
$f\dvtx\mathbb{R}_+\to\mathbb{R}_+$,
%
%
\begin{equation}\label{eqlimVq}V_{q}f(0):=\lim_{x\to
0+}V_{q}f(x)=\int^\infty_{0}\e^{R,H}_{0+}
(e^{-qR_{t}}H^{-\alpha}_{t}\kappa_{q}(H_{t},f))\,\dd t
\end{equation}
and, in particular,
%
%
\begin{equation}\label{unity2}
V_{q}1(0)=\frac{1}{q}.
\end{equation}
To finish this part of the proof we will now describe $V_{q}f(0)$.
Using the self-similarity of $(R,H)$, making a change of variables and
from the identity for the law of $(R_{1},H_{1})$ under $\p^{R,H}_{0+}$
obtained in Theorem~\ref{limitladder}, we obtain that
%
%
\begin{eqnarray}\label{eqlimVq2}
V_{q}f(0)&=&\int^\infty_{0} \e^{R,H}_{0+}(e^{-qtR_{1}}
(t^{1/\alpha}H_{1})^{-\alpha}\kappa_{q}(t^{1/\alpha
}H_{1},f))\,\dd t\nonumber\\
&=&\alpha\int^\infty_{0}\e^{R,H}_{0+}(\exp\{-qt^{\alpha
}R_{1}H^{-\alpha}_{1}\}H^{-\alpha}_{1}\kappa_{q}(t,f)
)\,\frac{\dd t}{t}\\
&=&\frac{\alpha}{\alpha\mu_{+}}\int^\infty_{0}\er(\exp
\{-qt^{\alpha}\widetilde{I}\}\kappa_{q}(t,f)
)\,\frac{\dd t}{t}.\nonumber
\end{eqnarray}
Next, we observe that the kernel $\kappa_{q}$ can be represented as
%
%
\begin{equation}\label{kusnetzov}
\kappa_{q}(z,f)=\int^\infty_{0}f(ze^{-x})(ze^{-x})^{\alpha}\er
^{\dagger}_{x}\biggl(\exp\biggl\{-qz^{\alpha}\int^{\varsigma
_0}_{0}e^{-\alpha\xi_{u}}\,\dd u\biggr\}\biggr)\widehat
{\mathcal{V}}(\dd x),\hspace*{-40pt}
\end{equation}
where we recall $\er^{\dagger}$ is the law of the process $\xi$
conditioned to stay positive, reflected at its future infimum, and
$\varsigma_0$ denotes its first hitting time of $0$. Indeed, this is
an easy consequence of the fact that the image under time reversal of
$\overline{n}$ is the excursion measure, say $n^{\uparrow}$, of the
process of excursions of $\xi^{\uparrow}$ from its future infimum
(see, e.g., Lemma 4 of~\cite{bertoindecomp}), that under $n^{\uparrow
}$ the coordinate process has the Markov property with semi-group
$P^{\dagger}_{t}f(x)=\er^{\dagger}_{x}(f(X_{t}), t<\varsigma
_{0})$, $t\geq0$, and the formula (\ref{dualrenewal}).

Hence, using identity (\ref{kusnetzov}), the former limit and making
some elementary manipulations we obtain that
%
%
\begin{eqnarray}\label{LTVq}
V_{q}f(0)
&=&\frac{\alpha}{\alpha\mu_{+}}\er\biggl(\int^{\infty
}_{0}e^{-qv^{\alpha}\widetilde{I}}\kappa_{q}(v,f)\,\frac{\mathrm
{d}v}{v}\biggr)\nonumber\\
&=&\frac{\alpha}{\alpha\mu_{+}}\er\biggl(\int^{\infty
}_{0}e^{-qv^{\alpha}\widetilde{I}}\int^\infty
_{0}f(ve^{-x})(ve^{-x})^{\alpha}\nonumber\\
&&\qquad\quad\hspace*{73pt}{}\times\er^{\dagger}_{x}\biggl(\exp\biggl\{
-qv^{\alpha}\int^{\varsigma_{0}}_{0}e^{-\alpha\xi_{u}}\,\dd u\biggr\}\biggr)\widehat
{\mathcal{V}}(\dd x)\,\frac{\dd v}{v}\biggr)\nonumber\\[-4pt]\\[-12pt]
&=&\frac{\alpha}{\alpha\mu_{+}}\int^{\infty}_{0} l^{\alpha
-1}f(l)\int^\infty_{0}\er(e^{-ql^{\alpha}e^{\alpha
x}\widetilde{I}})\nonumber\\
&&\qquad\quad\hspace*{70pt}{}\times\er^{\dagger}_{x}\biggl(\exp\biggl\{
-ql^{\alpha}e^{\alpha x}\int^{\varsigma_{0}}_{0}e^{-\alpha\xi
_{u}}\,\dd u\biggr\}\biggr)\widehat{\mathcal{V}}(\dd x)\,\dd l\nonumber\\
&=&\alpha\int^{\infty}_{0}
l^{\alpha-1}f(l)\biggl(\int^\infty_{0}e^{-ql^{\alpha}x}\eta
(\dd x)\biggr)\,\dd l,\nonumber
\end{eqnarray}
where $\eta$ is the sigma finite measure defined by
\[
\eta(f)=\frac{1}{\alpha\mu_{+}}\int_{\re^{3}_{+}}\pr^{\dagger
}_{x}\biggl(\int^{\varsigma_{0}}_{0}e^{-\alpha\xi_{u}}\,\dd u\in\dd
s\biggr)f\bigl(e^{\alpha x}(t+s
)\bigr)\pr(\widetilde{I}\in\dd t)\widehat{\mathcal
{V}}(\dd x).
\]
Note that on account of (\ref{unity2}) we also have that $\int
_{\mathbb{R}_+} x^{-1}\eta(\dd x)=1$. To complete the proof we
invert the Laplace transform in (\ref{LTVq}) and recover that the
entrance law of $X$ under $\p_{0+}$ is given by
\[
\e_{0+}(f(X_{t}))=\int_{\mathbb{R}_+}f(t^{1/\alpha
}x^{-1/\alpha})x^{-1}\eta(\dd x).
\]

\subsection{\texorpdfstring{Proof of the weak convergence in Theorem \protect\ref{thm:mainconvresult}}{Proof of the weak convergence in Theorem 4}}\label{weakconvergence}

In Theorems~\ref{limitladder} and~\ref{thm:mainconvresult}, we proved
the convergence in the sense of finite-dimensional distributions for
$(R,H)$ and $X$ under $\p_x$, as
$x\downarrow0$. As a consequence, we deduce the following corollary
which corresponds to the last part of the statement of Theorem~\ref
{thm:mainconvresult} and which was already obtained in~\cite{CCh}
under the
additional hypothesis $\mathbf{E}(\log^+\times\int_0^{\varsigma
_{(1,\infty)}}\exp\xi_s\, \dd s)<\infty$.
\begin{corollary}
Under the conditions of Theorem~\ref{thm:mainconvresult}, the family
of probability measures $(\p_x)$ converges weakly
toward $\p_{0+}$ as $x$ tends to $0$.
\end{corollary}
\begin{pf} Fix a sequence $(x_n)_{n\ge1}$ of positive real numbers which
converges to~$0$. Recall from Theorem~\ref{thm:mainconvresult} that
the sequence of probability measures $(\p_{x_n})$
converges to $\p_{0+}$ in the sense of finite-dimensional
distributions as $n\rightarrow\infty$. We will first prove that
the sequence $(\p_{x_n})$ actually converges weakly on $\mathbb{D}([0,1])$.

To this end, we apply Theorem 15.4 of~\cite{Bi}. First since $X$ has
$\p_{0+}$-a.s. no fixed discontinuities, the condition
$\p_{0+}(X_{1-}=X_1)=1$ is satisfied. Then for $0<\delta<1$, define
\[
W(\delta)=\sup\min\{|X_t-X_{t_1}|,|X_{t_2}-X_t|\} ,
\]
where the supremum extends over $t_1$, $t$ and $t_2$ in $[0,1]$ satisfying
\[
t_1\le t\le t_2 \quad\mbox{and}\quad t_2-t_1\le\delta.
\]
From Theorem 15.4 of~\cite{Bi}, it remains to prove that
for all $\varepsilon>0$ and $\chi>0$, there exist $0<\delta<1$ and
an integer $n_0$ such that for all $n\ge n_0$
%
%
\begin{equation}\label{second}
\p_{x_n}\bigl(W(\delta)>\varepsilon\bigr)\le\chi.
\end{equation}
Let $\gamma,\delta\in(0,1)$ such that $\delta<\gamma$ and note
that for
$t_1$, $t$ and $t_2$ in $[0,1]$ satisfying $t_1\le t\le t_2$ and
$t_2-t_1\le\delta$, if $t_1\in[0,\gamma)$,
then $\min\{|X_t-X_{t_1}|,|X_{t_2}-X_t|\}\le\sup_{0\le t\le2\gamma
}X_t$. Hence
\[
W(\delta)\le\sup_{0\le t\le2\gamma}X_t+\sup\min\{
|X_t-X_{t_1}|,|X_{t_2}-X_t|\} ,
\]
where the supremum extends over $t_1$, $t$ and $t_2$ in $[\gamma,1]$ satisfying
\mbox{$t_1\le t\le t_2$} and $t_2-t_1\le\delta$. The second term on the
right-hand side of the above inequality
is smaller than $W(\delta)\circ\theta_\gamma$, so that
%
%
\begin{equation}\label{eneq}
W(\delta)\le\sup_{0\le t\le2\gamma}X_t+W(\delta)\circ\theta
_\gamma.
\end{equation}
From (\ref{eneq}) and the Markov property one has for all $n\ge1$
and for all $\delta,\gamma$ as above,
%
%
\begin{equation}\label{eneq2}\qquad
\p_{x_n}\bigl(W(\delta)>\varepsilon\bigr)\le\p_{x_n}\Bigl(\sup_{0\leq t\le
2\gamma}X_t>\varepsilon/2\Bigr)+\e_{x_n}\bigl(\p_{X_\gamma
}\bigl(W(\delta)>\varepsilon/2\bigr)\bigr) .
\end{equation}
To deal with the first term in (\ref{eneq2}), we pick $u>0$ and we write
\begin{eqnarray*}
\p_{x_n}\Bigl(\sup_{0\leq t\le2\gamma}X_t>\varepsilon/2
\Bigr)&=&\p_{x_n}\Bigl(\sup_{0\le t\leq2\gamma}X_t>\varepsilon/2,
R_u<2\gamma\Bigr)\\
&&{}+\p_{x_n}\Bigl(\sup_{0\leq t\le2\gamma}X_t>\varepsilon/2, R_u\ge
2\gamma\Bigr)\\
&\le&\p_{x_n}(R_u<2\gamma)+\p_{x_n}
(H_u>\varepsilon/2) .
\end{eqnarray*}
From Theorem~\ref{limitladder}, $\p_{x_n}(R_u\in\dd s, H_u\in
\dd y)$ converges weakly
to $\p_{0+}(R_u\in\dd s$, $H_u\in\dd y)$, as $n$ tends to
$\infty$, hence $\lim_n\p_{x_n}(R_u< 2\gamma)=\p_{0+}(R_u< 2\gamma)$
and $\lim_n\p_{x_n}(H_u>\varepsilon/2)=\p_{0+}(H_u>\varepsilon
/2)$ (without loss of generality we can make sure that
$2\gamma$ and $\varepsilon/2$ are points of continuity
of the distribution functions of $H_u$ and~$R_u$, resp., under $\p
_{0+}$). Moreover since $H_0=0$, $\p_{0+}$
a.s., we have $\lim_{u\downarrow0}\p_{0+}(H_u>\varepsilon/2)=0$, so
we may find $u>0$ and an integer $n_1$ such that for all $n\ge n_1$,
$\p_{x_n}(H_u>\varepsilon/2)<\chi/4$. Then since $R_u>0$, $\p
_{0+}$-a.s., we may find $\gamma\in(0,1)$ and an integer $n_2$ such
that for all $n\ge n_2$, $\p_{x_n}(R_u<2\gamma)\le\chi/4$.

Next we deal with the second term in (\ref{eneq2}). From Theorem \ref
{thm:mainconvresult}, $\p_{x_n}(X_{\gamma}\in dz)$
converges weakly toward
$\p_{0+}(X_{\gamma}\in dz)$ as $n$ tends to $\infty$. Moreover, we
may easily check, using the Lamperti representation and general properties
of L\'{e}vy processes, that $x\mapsto\p_x(W(\delta)>\varepsilon/2)$
is continuous on $(0,\infty)$. Since $\p_{0+}(X_{\gamma}=0)=0$, we have
$\lim_{n\rightarrow\infty}\e_{x_n}(\p_{X_\gamma}(W(\delta
)>\varepsilon/2))=\e_{0+}(\p_{X_\gamma}(W(\delta)>\varepsilon/2))$.
For all $x\ge0$, $\lim_{\delta\downarrow0}W(\delta)=0$, $\p
_x$-a.s. (see pages 110 and 119 of~\cite{Bi}) so that using
dominated convergence, $\lim_{\delta\downarrow0}\e_{0+}(\p
_{X_\gamma}(W(\delta)>\varepsilon/2))=0$.
Then, we may find $n_3$ and $\delta$ such that for all $n\ge n_3$, $\e
_{x_n}(\p_{X_\gamma}(W(\delta)>\varepsilon/2))\le\chi/2$.

We conclude that (\ref{second}) is satisfied with $\delta$ and
$n_0=\max(n_1,n_2,n_3)$, so that the sequence $(\p_{x_n})$
restricted to $\mathbb{D}([0,1])$ converges weakly to $\p_{0+}$. Then
it follows from the same arguments that
the sequence $(\p_{x_n})$ restricted to $\mathbb{D}([0,t])$ converges
weakly to $\p_{0+}$, for each $t>0$. Finally it
remains to apply Theorem 16.7 of~\cite{Bi} to conclude that $(\p
_{x_n})$ converges weakly to
$\p_{0+}$ on $\mathbb{D}([0,\infty))$.
\end{pf}


%

%
\printaddresses

\end{document}